\documentclass[a4paper, 12pt]{article}
\usepackage{latexsym}
\usepackage[T1]{fontenc}
\usepackage{amsmath}
\usepackage{amssymb}
\usepackage{graphicx,graphpap}
\usepackage{color,xcolor}
\usepackage{kbordermatrix}
\usepackage{dsfont}
\newtheorem{defn}{Definition}[section]
\newtheorem{lem}[defn]{Lemma}
\newtheorem{prop}[defn]{Proposition}

\newtheorem{cor}[defn]{Corollary}
\newtheorem{rem}[defn]{Remark}
\newcommand{\vc}[1]{\boldsymbol{#1}}

\newcommand{\cqfd}{\hfill $\square$}

\begin{document}

\title{The time-dependent expected reward and deviation matrix of a finite QBD process}
\author{Sarah Dendievel\thanks{Ghent University, Department of
  Telecommunications and Information Processing, SMACS Research Group,
  Sint-Pietersnieuwstraat 41, B-9000 Gent, Belgium, \texttt{Sarah.Dendievel@UGent.be}}, Sophie Hautphenne\thanks{Ecole Polythechnique F\'ed\'erale de Lausanne,  Institute of Mathematics, and The University of Melbourne, School of Mathematics and Statistics, \texttt{sophiemh@unimelb.edu.au} }, Guy Latouche\thanks{Universit\'e libre de Bruxelles,  Facult\'e des sciences, CP212, Boulevard du Triomphe 2, 1050 Bruxelles,
  Belgium,        \texttt{latouche@ulb.ac.be}}, \\and Peter Taylor\thanks{The University of Melbourne, School of Mathematics and Statistics,  \texttt{pgt@ms.unimelb.edu.au} }}
\maketitle

\begin{abstract}
Deriving the time-dependent expected reward function associated with a continuous-time Markov chain involves the computation of its transient deviation matrix. In this paper we focus on the special case of a finite quasi-birth-and-death (QBD) process, motivated by the desire to compute the expected revenue lost in a MAP/PH/1/C queue. 

We use two different approaches in this context. The first is based on the solution of a finite system of matrix difference equations; it provides an expression for the blocks of the expected reward vector, the deviation matrix, and the mean first passage time matrix. The second approach, based on some results in the perturbation theory of Markov chains, leads to a recursive method to compute the full deviation matrix of a finite QBD process. We compare the two approaches using some numerical examples. 
\end{abstract}

\section{Introduction}Our analysis of the expected reward in a finite QBD process is motivated by the problem
of computing the expected amount of lost revenue in a \emph{MAP/PH/1/C} queue over a finite time horizon $[0,t]$, given its initial occupancy. Since MAPs and PH distributions are the most general matrix extensions of Poisson processes and exponential distributions, respectively,  we can think of this problem as a matrix generalisation of the similar analysis for the $M/M/C/C$ model considered in Chiera and Taylor \cite{chiera2002unit} and the $M/M/1/C$ model in Braunsteins, Hautphenne and Taylor \cite{braunsteins2016}. In these queueing models, customers are lost when they arrive to find $C$ customers already present. Assuming that each arriving customer brings a certain amount of revenue, we are interested in calculating the expected amount of revenue that the queue will lose over a finite time horizon $[0,t]$, as well as exploring the limit of the rate of losing revenue in the asymptotic regime. 

Solving the expected lost revenue problem is important, for example, if we wish to find a way of managing a system where a number of, possibly different, queues share {a number of} servers. If it is feasible to reallocate servers from one queue to another every $t$ time units, then a rational method for performing the allocation is for each queue to observe its occupancy at time 0, calculate the expected revenue lost in time $[0,t]$ for a range of capacities, given its initial occupancy, and then to allocate the servers to minimise the total expected amount of lost revenue over $[0,t]$. At time $t$, the calculation can be performed again, based upon the occupancies at that time and a reallocation performed if it is optimal to do so.

The \emph{MAP/PH/1/C} queue can be modelled as a finite QBD process with generator matrix $Q$, and levels $0,1,\ldots,C$ corresponding to the possible queue lengths. If $\vc R(t)$ is a vector containing the expected revenue lost in $[0,t]$ conditional on the initial state (level and phase) of the system, then computing $\vc R(t)$ reduces to solving a special case of the time-dependent version of Poisson's equation of the form 
\begin{eqnarray}\nonumber  \vc R(0)&=&\vc 0\\ \label{eq0}\vc R'(t)&=& Q \vc R(t)+\vc g,\end{eqnarray} where $\vc 0$ is a column vector of $0$'s, and $\vc g$ is a column vector containing the reward (loss) per unit of time in each state of the system. Since the \emph{MAP/PH/1/C} system loses revenue only when it is at full capacity, the only non-zero entries of $\vc g$ in our motivating example are those corresponding to level $C$. 

The solution of \eqref{eq0}, given by
\begin{equation}\label{eq3}\vc R(t)=(\vc\pi\vc g)\vc 1\,t+D(t)\,\vc g,\end{equation} 
where $\vc\pi$ is the stationary vector of the QBD process and $\vc 1$ denotes the column vector of 1's, involves the transient deviation matrix, 
\begin{equation}
\label{eq:dt}
{D(t)}=\int_0^t \left(e^{Q u}-\vc1\vc\pi\right) du,
\end{equation}
see \cite{braunsteins2016}. As $t\rightarrow\infty$, $D(t)$ converges to the deviation matrix $D$ discussed in Coolen-Schrijner and van Doorn \cite{CoolenSc}, which corresponds to the group inverse of $-Q$, and the expected lost revenue function has a linear asymptote, $\vc R(t)\sim(\vc\pi\vc g)\vc 1\,t+D\vc g.$ 

After providing more detail on the reward function $\vc R(t)$ associated with a QBD process in Section 2, we tackle the computation of $\vc R(t)$ in transient and asymptotic regimes, and the corresponding matrices $D(t)$ and $D$, using two different approaches, each having some advantages in comparison to the other. In the first approach, developed in Section 3, we assume that the QBD process is non null-recurrent, and we place ourselves in the general context where the reward vector $\vc g$ is not restricted to any particular structure.
We use systems of matrix difference equations to gain insight into the block-structure of the vector $\vc R(t)$ and of the matrices $D(t)$ and $D$. 
In addition, this method also provides us with the blocks of the matrix of the mean first passage times in the QBD process. We obtain simple expressions for the relevant quantities that highlight the role of the maximal capacity $C$. We also derive limiting results when the maximal capacity increases to infinity. This approach is effective if one wants to focus on particular blocks of $\vc R(t)$, $D(t)$ and $D$ rather than on the full matrices, as is the case in our motivating example. In practice, it also avoids dealing with large matrices: the size of the matrices involved in the expressions is at most twice the size of the phase space of the QBD process, regardless the value of $C$.

 The second approach, described in Section 4, relies on some elegant results from the perturbation theory of Markov chains, and leads to a recursive formula expressing the full deviation matrices $D(t)$ and $D$ for a given value of $C$ in terms of the corresponding matrices for a system with capacity $C-1$. 
 
 The two approaches complement each other and their use may depend on the context: the first method is algebraic because we solve matrix difference equations, but the solutions have a probabilistic interpretation, while the second method is probabilistic but leads to a solution with an algebraic flavour.


In Section 5, we provide numerical illustrations, starting with our motivating \emph{MAP/PH/1/C} example which loses revenue only at full capacity, and moving on to look at a case where the reward function $\vc g$ is non-zero for all levels of the system. We compare the computational complexity of the two approaches and show that there is always a threshold value of $C$ at which one method surpasses the other in terms of CPU time.



\section{Background}

Let $\{X(t):t\geq 0\}$ be an ergodic continuous-time Markov chain on a finite state-space, with generator $Q$ and stationary distribution $\vc\pi$. 
The \emph{deviation matrix} of $\{X(t)\}$ is the matrix
\begin{equation}
\label{eq:d}
{D}=\int_0^\infty \left(e^{Q u}-\vc1\vc\pi\right) du,
\end{equation}
whose components may be written as
${D}_{ij}=\lim_{t\rightarrow\infty}[N_{ij}(t)-N_{\boldsymbol{\pi}j}(t)],$
where $N_{ij}(t)$ is the expected time spent in state $j$ during
the interval of time $[0,t]$ given that the initial state is $i$,
and $N_{\boldsymbol{\pi}j}(t)$ is the same quantity but conditional on the initial state
having the distribution $\vc\pi$ (see Da Silva Soares and Latouche \cite{da2002group}). 

The \emph{group inverse}  $A^{\#}$ of a matrix $A$, if it exists, is defined as the unique solution to $A A^{\#} A=A$, $A^{\#}A A^{\#}=A^{\#}$, and $A^{\#} A=A A^{\#}$. 
From Campbell and Meyer \cite[Theorem 8.5.5]{campbell2009generalized}, the group inverse $Q^{\#}$ of the infinitesimal generator $Q$ of any finite
Markov chain  is the unique solution of the two equations 
\begin{align}
Q Q^{\#} & =I-W,\label{Cesaro}\\
W Q^{\#} & =0.\label{Cesaro2}
\end{align}
where $W=\lim_{t\rightarrow\infty}[\exp(Qt)]$. When it exists, the deviation matrix is related to the group inverse $Q^{\#}$ of $Q$ by the relation $$D=-Q^{\#}.$$ 
If $Q$ is irreducible,
then $D$ is also the unique solution of the system 
\begin{align}
QD & =\boldsymbol{1\pi}-I,\label{eq:deveq2}\\
\boldsymbol{\pi}D & =\vc0;\label{eq:deveq1}
\end{align}
see \cite{CoolenSc} for more detail.
In addition, $D$ also satisfies 
\begin{equation}
D\boldsymbol{1}=\vc0.\label{EqofQsharpe}
\end{equation}


We see that $D=\lim_{t\rightarrow\infty}D(t)$, where $D(t)$ is the transient deviation matrix defined by (\ref{eq:dt}). Properties of the transient deviation matrix are discussed in \cite{braunsteins2016}. 

Let $V_{ij}(t)$ be the expected cumulative time spent in state $j$  in the time interval $[0,t]$, given that the process starts in state $i$ at time $0$, and define the matrix $V(t)=(V_{ij}(t))$.  It follows that
 $$V(t)=\int_0^t e^{Qu}du=\vc1\vc\pi \,t+{D(t)},$$
and $V(t)$ has the linear asymptote $\bar{V}(t)=\vc1\vc\pi \,t+D.$
If we associate a reward (or loss) $g_j$ per time unit when the Markov chain $\{X(t)\}$ occupies state $j$, and define the vector $\vc g:=(g_j)$, then the expected cumulative reward up to time $t$, given that the chain starts in state $i$, is given by
$R_i(t):=(V(t)\vc g)_i,$ so that the vector $\vc R(t):=(R_i)$ satisfies
\begin{equation}\label{eq3}\vc R(t)=(\vc\pi\vc g)\vc1\, t+{D(t)}\vc g,\end{equation} and $\vc R(t)$ has the linear asymptote \begin{equation}\label{lin_ass}\bar{\vc R}(t)=(\vc\pi\vc g)\vc1\, t+{D}\vc g.\end{equation} Observe that \eqref{eq3} is the solution of a finite horizon version of Poisson's equation,
\begin{equation}\label{eq1}\vc R'(t)= Q \vc R(t)+\vc g \end{equation}with $\vc R(0)=\vc 0$.  In the Laplace transform domain where, for $\Re(s)>0$, $$\tilde{\vc R}(s)=\int_0^\infty e^{-st} \vc R(t)\, dt,$$\eqref{eq3} and \eqref{eq1} become respectively
\begin{eqnarray}\label{LTTD}\tilde{\vc R}(s)&=&(1/s^2)\, (\vc\pi\vc g) \vc 1+\tilde{D}(s)\,\vc g\\\label{sol}
 &=& (sI-Q)^{-1} (1/s)\vc g,
\end{eqnarray}where $\tilde{D}(s)$ is the Laplace transform of the transient deviation matrix, given by
\begin{equation}\label{ltd}\tilde{D}(s)= (1/s) (sI-Q)^{-1}-(1/s^2)\vc 1\vc\pi. \end{equation} {Another equivalent expression, obtained after some algebraic manipulations, is given by
$$\tilde{D}(s)= (1/s) (sI-Q)^{-1}(I-\vc 1\vc\pi). $$}

In this paper, we are interested in computing the expected reward function $\vc R(t)$, and the associated deviation matrices $D(t)$ and $D$, for any $t\geq 0$, when the Markov chain corresponds to a finite (level-independent) QBD process. Such a process is a two-dimensional continuous-time Markov chain $\{\vc X(t)=(J(t),\varphi(t)), t\geq 0\}$, where the variable $J(t)$, taking values in $\{0,1,\ldots,C\}$, is called the \textit{level} of the process at time $t$, and the variable $\varphi(t)$ taking values in $\{1,2,\ldots,n\}$ is called the \textit{phase} of the process at time $t$. The generator of the QBD process has a block-tridiagonal form given by
\begin{equation}\label{GeneratorQ}
Q=\kbordermatrix{\mbox{}&0&1&2&3&\ldots& C-1&C\\
0&B_0&A_1&& &&&\\
1&A_{-1}&A_0&A_1& &&&\\
2&&A_{-1}&A_0& A_1&&&\\
\vdots& & && &\ddots&&\\C-1&&&& &A_{-1}&A_{0}&A_1\\
C&&&& &&A_{-1}&C_0\\
},
\end{equation}where the block matrices are of size $n\times n$.

The vectors $\vc R(t)$ and $\vc g$ corresponding to the QBD process are of length $(C+1)n$ and can be structured into block sub-vectors $\vc R_k(t)$ and $\vc g_k$ of length $n$, for 
$0\leq k\leq C$, corresponding to each level of the QBD process. The same block decomposition holds for the deviation matrices $D(t)$ and $D$, and we write $D_{k,\ell}(t)=\left(D_{(k,i)(\ell,j)}(t)\right)_{1\leq i,j\leq n}$ ($D_{k,\ell}=\left(D_{(k,i)(\ell,j)}\right)$ respectively) for the $(k,\ell)$th block of these matrices.



The next two sections address the following questions:
\begin{itemize}
\item[\textit{(i)}] How do we use the structure of the QBD process to compute the blocks of ${\vc R}(t)$, $D(t)$ and $D$?
\item[\textit{(ii)}] How do the deviation matrices $D(t)$ and $D$ differ for two successive values of the capacity $C$? \label{questiondeux}
\end{itemize}

\section{A matrix difference equation approach}\label{mda}In this section, we answer Question \emph{(i)}.
We deal with the transient regime first, followed by the asymptotic regime.

\subsection{Transient regime}\label{tr}We analyse the transient behaviour of the expected reward function in the Laplace transform domain.
An explicit expression for $\tilde{\vc R}(s)$ is given in \eqref{sol} but 
this expression does not give any insight in the structure of the blocks $\tilde{\vc R}_k(s)$, which are conditional on the initial queue length, nor in the role of the maximum capacity $C$ in the solution. Instead, we rewrite \eqref{sol} as a system of second-order matrix difference equations
\begin{eqnarray}\label{Rt0}
(B_0-sI)\tilde{\vc R}_0(s)+ A_1\,\tilde{\vc R}_1(s)&=&- \vc g_0(s)\\\label{Rtn}
A_{-1}\,\tilde{\vc R}_{k-1}(s)+(A_0-sI)\tilde{\vc R}_k(s)+ A_1\,\tilde{\vc R}_{k+1}(s)&=&- \vc g_k(s),\; 1\leq k\leq C-1\;\qquad\\\label{RtC}A_{-1}\,\tilde{\vc R}_{C-1}(s)+(C_0-sI)\tilde{\vc R}_C(s)&=&- \vc g_C(s),
\end{eqnarray}where $\vc g_k(s)=(1/s) \,\vc g_k$, and we focus on the computation of the blocks $\tilde{\vc R}_k(s)$, for $0\leq k\leq C$.

  For any $s\geq 0$, let 
$G(s)$ be the minimal nonnegative solution to the matrix quadratic equation
\begin{equation}
\label{eqG}A_{-1}+(A_0-sI)X+A_1X^2=0,\end{equation} and let $\hat{G}(s)$ be the minimal nonnegative solution to
\begin{equation}\label{eqGh}A_{1}+(A_0-sI)X+A_{-1}X^2=0.\end{equation} The matrices $G(s)$ and $\hat{G}(s)$ can be computed numerically using any of the linear or quadratic algorithms discussed in Latouche and Ramaswami \cite{latram}. Assuming that the process is not restricted by an upper boundary, for $s>0$, the $(i,j)th$ entry of the matrix $G(s)$ contains the Laplace transform 
\begin{equation}
\label{eq:gs}
\mathbb{E}\left[e^{s \theta_{k-1}}\,\mathds{1}\{\varphi(\theta_{k-1})=j\}\,|\,J(0)=k, \varphi(0)=i\right]
\end{equation}
where $\theta_{k-1} = \inf_{t>0} \{J(t)=k-1\}$
and $\mathds{1}\{\cdot\}$ denotes the indicator function.
Similarly, assuming that the process is not restricted by a lower boundary, the $(i,j)th$ entry of the matrix $\hat G(s)$ contains the Laplace transform 
\begin{equation}
\label{eq:hatgs}
\mathbb{E}\left[e^{s \theta_{k+1}}\,\mathds{1}\{\varphi(\theta_{k+1})=j\}\,|\,J(0)=k, \varphi(0)=i\right].
\end{equation}
Both $G(s)$ and $\hat{G}(s)$ are sub-stochastic, and therefore have a spectral radius strictly less than 1, as long as $s>0$. The first passage probability matrices $G$ and $\hat{G}$ correspond to $G(0)$ and $\hat{G}(0)$. When the unrestricted QBD process is transient, $G$ is sub-stochastic and $\hat{G}$ is stochastic, and conversely, when the unrestricted QBD process is positive recurrent, $G$ is stochastic and $\hat{G}$ is sub-stochastic. Finally, we let $$H_0(s)=-(A_0-sI+A_1G(s)+A_{-1}\hat{G}(s))^{-1},$$ which is well defined for any $s>0$, and
\begin{equation}\label{ho}H_0=-(A_0+A_1G+A_{-1}\hat{G})^{-1},\end{equation} which is well defined in the non null-recurrent case. 

{In the sequel, we use the convention that an empty sum, such as $\sum_{j=1}^0 $ or $\sum_{j=0}^{-1} $, is zero.}


\begin{lem} \label{solhomu}For any $C\geq 1$ and $1\leq k\leq  C$, the general solution of the second-order matrix difference equation
\begin{equation}\label{equu}A_{-1}\,\vc u_{k-1}+A_0\,\vc u_k+A_1\,\vc u_{k+1}=-\vc g_k\end{equation} is given by
\begin{equation}
\label{gensoluk}\vc u_k=G^k\,\vc v+\,\hat{G}^{C-k}\vc w+\vc\nu_k(C),\end{equation} where $\vc v$ and $\vc w$ are arbitrary vectors, and 
$$\vc\nu_k(C)=\sum_{j=0}^{k-1} G^j H_0 \,\vc g_{k-j}+\sum_{j=1}^{C-k} \hat{G}^j H_0\,\vc g_{k+j}.$$
\end{lem}

\noindent \textit{Proof.} First, using \eqref{eqG}, \eqref{eqGh} and \eqref{ho} we can show that \eqref{gensoluk} is a solution to \eqref{equu} for any arbitrary vectors $\vc v,\vc w\in\mathbb{C}^{n}$. It remains to show that all solutions of \eqref{equu} can be expressed in the form of \eqref{gensoluk}.

Observe that there exists a nonsingular matrix $M$ such that 
$$\hat{G}M=MJ$$with $$J=\left[\begin{array}{cc} V & 0\\ 0& W\end{array}\right],$$ where $V$ is a non-singular square matrix of order $p$, and $W$ is a square matrix of order $q$ with $p+q=n$ and sp$(W)=0$. For instance, we may choose $M$ such that $J$ is the Jordan normal form of $\hat{G}$; in this case, $V$ contains all the blocks for the non-zero eigenvalues of $\hat{G}$ (which lie within the unit circle), and $W$ contain all the blocks for the zero eigenvalues. The matrices $M$ and $M^{-1}$ may be partitioned as $$M=\left[\begin{array}{c|c}L&K\end{array}\right],\qquad M^{-1}=\left[\begin{array}{c}E\\\hline F\end{array}\right],$$ where $L$ has dimension $n\times p$, $E$ has dimension $p\times n$, and $EL=I$, $FL=0$.
 Lemma 10 in Bini \textit{et al.} \cite{bini} states that the general solution of  \eqref{equu} is given, for $k\geq 1$, by
\begin{equation}\label{solukbini}\vc u_k=G^k\,\vc v+\,LV^{-k}\,\vc z+\vc\sigma_k,\end{equation} where $\vc v\in\mathbb{C}^{n}$ and $\vc z\in\mathbb{C}^p$ are arbitrary vectors, and 
\begin{equation}\label{sigk0}\vc\sigma_k=\sum_{j=0}^{k-1}(G^j-LV^{-j}E)H_0\,\vc g_{k-j}+\vc\tau_k,\end{equation} with $$\vc\tau_k=\sum_{j=1}^{\nu-1} K W^j F H_0 \,\vc g_{k+j},$$and $\nu$ the smallest integer such that $M^\nu=0$ ($\vc\tau_k=\vc 0$ if $\hat{G}$ is invertible). 

Since the eigenvalues of $V$ lie within the unit circle, the negative powers of $V$ appearing in \eqref{solukbini} and \eqref{sigk0} have entries that can take unbounded values, which can lead to numerical instabilities. Next, we rewrite \eqref{solukbini} in a more convenient form in order to get rid of any negative power of $V$. First observe that \eqref{solukbini} is equivalent to
\begin{eqnarray*}\label{solukbini2}\vc u_k&=&G^k\,\vc v+\sum_{j=0}^{k-1}G^jH_0\,\vc g_{k-j}+\,LV^{-k}\,(\vc z-\sum_{j=0}^{k-1}V^{k-j}EH_0\,\vc g_{k-j})+\vc\tau_k,\end{eqnarray*}where the negative powers of $V$ appear in the term  
\begin{eqnarray*}\vc s_k&:=&LV^{-k}\,(\vc z-\sum_{j=0}^{k-1}V^{k-j}EH_0\,\vc g_{k-j})=LV^{-k}\,(\vc z-\sum_{j=1}^{k}V^{j}EH_0\,\vc g_{j}).\end{eqnarray*} Letting $\vc z=\vc y +\sum_{j=1}^{C}V^{j}EH_0\,\vc g_{j}$, where $\vc y\in\mathbb{C}^p$ is an arbitrary vector, we obtain
\begin{eqnarray*}\vc s_k&=
LV^{-k}\,\vc y + \sum_{j=1}^{C-k}LV^{j}EH_0\,\vc g_{k+j}.
\end{eqnarray*}
 Let us fix a vector $\vc y\in \mathbb{C}^p$. Then, for any $k\leq C$, $LV^{-k}\,\vc y$ can be equivalently written as $\hat{G}^{C-k} \vc w$ with $\vc w=LV^{-C} \vc y$. Indeed, 
 \begin{eqnarray*}\hat{G}^{C-k} \vc w&=&\hat{G}^{C-k} LV^{-C} \vc y\\&=&[L V^{C-k} E+K W^{C-k}F]LV^{-C}\vc y\\&=&L V^{C-k}V^{-C}\vc y\\&=&L V^{-k}\vc y.
 \end{eqnarray*}Therefore, we have
$$\vc s_k= \hat{G}^{C-k} \vc w+\sum_{j=1}^{C-k}LV^{j}EH_0\,\vc g_{j+k},$$ and the general solution \eqref{solukbini} takes the form $$\vc u_k=G^k\,\vc v+\hat{G}^{C-k} \vc w+\vc \nu_k(C)$$ where \begin{eqnarray*} \vc \nu_k(C)&=&\sum_{j=0}^{k-1}G^jH_0\,\vc g_{k-j}+\sum_{j=1}^{C-k}LV^{j}EH_0\,\vc g_{k+j}+\sum_{j=1}^{\nu-1} K W^j F H_0 \,\vc g_{k+j}.\end{eqnarray*}Finally, note that $\vc g_{\ell}$ is defined for $0\leq \ell\leq C$ only, so we can set $\vc g_{k+j}=\vc 0$ for any $j>C-k$. With this, we have \begin{eqnarray*}\sum_{j=1}^{C-k}LV^{j}EH_0\,\vc g_{k+j}+\sum_{j=1}^{\nu-1} K W^j F H_0 \,\vc g_{k+j}&=&\sum_{j=1}^{C-k}\hat{G}^{j}H_0\,\vc g_{k+j},\end{eqnarray*}which shows that any solution to \eqref{equu} can be written in the form \eqref{gensoluk}.
 \hfill\cqfd
 \medskip
 
 The advantage of the solution \eqref{gensoluk} over the solution \eqref{solukbini} from \cite{bini} is that it does not require any spectral decomposition of the matrix $\hat{G}$, nor any matrix inversion (as long as $C\geq k$). Since the spectral radii of $G$ and $\hat{G}$ are bounded by 1, all matrix powers involved in \eqref{gensoluk} are bounded, and the computation of the solution is therefore numerically stable.  
 
 Lemma \ref{solhomu} naturally extends to the Laplace transform domain, as stated in the next Corollary.

\begin{cor}\label{gensolR} For any $s>0$, $C\geq 1$, and $1\leq k\leq C$, the general solution of the second-order matrix difference equation
\begin{equation}\label{Rtn2}A_{-1}\,{\vc u}_{k-1}(s)+(A_0-sI)\,{\vc u}_k(s)+A_1\,{\vc u}_{k+1}(s)=-\vc g_k(s)\end{equation} is given by
\begin{equation}\label{gen}{\vc u}_k(s)=G(s)^k\,\vc v(s)+\hat{G}(s)^{C-k}\,\vc w(s)+\vc\nu_k(s,C),\end{equation} where $\vc v(s)$ and $\vc w(s)$ are arbitrary vectors, and 
\begin{equation}\label{nuk}
\vc\nu_k(s,C)=\sum_{j=0}^{k-1} G(s)^j H_0(s) \,\vc g_{k-j}(s)+\sum_{j=1}^{C-k} \hat{G}(s)^j H_0(s)\,\vc g_{k+j}(s).\end{equation}
\hfill\cqfd
\end{cor}


We use Corollary \ref{gensolR} to obtain a closed-form expression for $\tilde{\vc R}_k(s)$ in terms of the matrices $G(s)$ and $\hat{G}(s)$ of the QBD process. 
\begin{prop}\label{prop26}
For any $s>0$, $C\geq 1$, and $0\leq k\leq C$, the Laplace transform of the expected reward function, conditional on the initial level $k$, is given by
\begin{equation}\label{sol2}\tilde{\vc R}_k(s)=G(s)^k \vc v(s,C) +\hat{G}(s)^{C-k}\vc w(s,C)+\vc\nu_k(s,C),\end{equation} where $\vc\nu_k(s,C)$ is given by \eqref{nuk}, and
$$\left[\begin{array}{c} \vc v(s,C)\\\vc w(s,C)\end{array}\right]=(-Z(s,C))^{-1}\left[\begin{array}{c}\vc g_0(s)+(B_0-sI)\,\vc\nu_0(s)+A_1\,\vc\nu_1(s)\\\vc g_C(s)+A_{-1}\,\vc\nu_{C-1}(s)+(C_0-sI)\,\vc\nu_C(s)\end{array}\right]$$with
\begin{equation}\label{ZsC}\small{
Z(s,C)=\left[\begin{array}{cc} (B_0-sI)+A_1G(s) & ((B_0-sI)\hat{G}(s)+A_1)\hat{G}(s)^{C-1}
\\(A_{-1}+(C_0-sI)G(s)) G(s)^{C-1}& A_{-1}\hat{G}(s)+(C_0-sI)
\end{array}\right].}\qquad\end{equation}In addition, $Z(s,C)$ can be written as$$Z(s,C)= \mathring{Q}(s,C)\left[\begin{array}{cc} I& \hat{G}(s)^{C}
\\G(s)^{C}&I
\end{array}\right],$$ where $\mathring{Q}(s,C)$ is the generator of the transient Markov chain obtained from the QBD process in which absorption can happen from any state at rate $s>0$, and restricted to levels 0 and $C$.
\end{prop}

\noindent\textit{Proof.} We apply Corollary \ref{gensolR} to \eqref{Rtn} and obtain that
the general solution of the second-order difference equation is given by
\eqref{sol2} for $1\leq k\leq C-1$. We then specify the arbitrary vectors $\vc v(s,C)$ and $\vc w(s,C)$ using the boundary conditions \eqref{Rt0} and \eqref{RtC}. Injecting the general solution into \eqref{Rt0} and \eqref{RtC} leads to the system of equations
\begin{eqnarray*}-\vc g _0(s)&=&(B_0-sI)(\vc v(s,C)+\hat{G}(s)^{C}\vc w(s,C)+\vc\nu_0(s))\\&&+A_1(G(s)\vc v(s,C)+\hat{G}(s)^{C-1}\vc w(s,C)+\vc\nu_1(s))\\-\vc g _C(s)&=&A_{-1}(G(s)^{C-1}\vc v(s,C)+\hat{G}(s)\vc w(s,C)+\vc\nu_{C-1}(s))\\&&+(C_0-sI)(G(s)^C\vc v(s,C)+\vc w(s,C)+\vc\nu_C(s)),
\end{eqnarray*} which can be rewritten as 
$$Z(s,C) \left[\begin{array}{c} \vc v(s,C)\\\vc w(s,C)\end{array}\right]=-\left[\begin{array}{c}\vc g_0(s)+(B_0-sI)\,\vc\nu_0(s)+A_1\,\vc\nu_1(s)\\\vc g_C(s)+A_{-1}\,\vc\nu_{C-1}(s)+(C_0-sI)\,\vc\nu_C(s)\end{array}\right],$$where
$Z(s,C)$ is given by \eqref{ZsC}.
To prove that the matrix $Z(s,C)$ is invertible, we now show that it can be written as the matrix product
\begin{equation}\label{e0}Z(s,C)=\mathring{Q}(s,C)\left[\begin{array}{cc} I& \hat{G}(s)^{C}
\\G(s)^{C}&I
\end{array}\right],\end{equation}where $\mathring{Q}(s,C)$ is the generator of the transient Markov chain obtained from the QBD process in which absorption can happen from any state at rate $s>0$, restricted to levels 0 and $C$, and is therefore invertible, and the matrix inverse 
$$\left[\begin{array}{cc} I& \hat{G}(s)^{C}
\\G(s)^{C}&I
\end{array}\right]^{-1}
$$
 exists because sp$(G(s))<1$ and sp$(\hat{G}(s))<1$  for any $s>0$.

The generator $\mathring{Q}(s,C)$ of the transient Markov chain in which absorption can happen from any state at rate $s>0$, restricted to levels 0 and $C$, can be written as $\mathring{Q}(s,C)=\mathcal{A}(s)+\mathcal{B} (-\mathcal{C}(s))^{-1}\mathcal{D}$ with
$$\mathcal{A}(s)=\left[\begin{array}{cc} (B_0-sI) & 0\\0& (C_0-sI)\end{array}\right],\quad \mathcal{B}=\left[\begin{array}{cccc} A_1 & 0&\ldots&0\\0&\ldots &0& A_{-1}\end{array}\right],$$
$$\mathcal{C}(s)=\left[\begin{array}{cccc}(A_0-sI) & A_1&&\\A_{-1}&\ddots&&\\&&\ddots&A_1\\&&A_{-1}&(A_0-sI)\end{array}\right],\quad \mathcal{D}=\left[\begin{array}{cc} A_{-1} & 0\\0&\vdots \\\vdots&0\\0& A_1\end{array}\right].$$ We then have 
\begin{eqnarray*}\lefteqn{\mathring{Q}(s,C)\left[\begin{array}{cc} I& \hat{G}(s)^{C}
\\G(s)^{C}&I
\end{array}\right]}\\&=&\mathcal{A}(s)\left[\begin{array}{cc} I& \hat{G}(s)^{C}
\\G(s)^{C}&I
\end{array}\right]+\mathcal{B} (-\mathcal{C}(s))^{-1} \mathcal{D} \left[\begin{array}{cc} I& \hat{G}(s)^{C}
\\G(s)^{C}&I
\end{array}\right],\end{eqnarray*} and to show that this is equal to $Z(s,C)$  given by \eqref{ZsC}, it suffices to show that 
\begin{equation}\label{toshow}\mathcal{B} (-\mathcal{C}(s))^{-1} \mathcal{D} \left[\begin{array}{cc} I& \hat{G}(s)^{C}
\\G(s)^{C}&I
\end{array}\right]= \left[\begin{array}{cc} A_1 G(s)& A_1\hat{G}(s)^{C-1}
\\A_{-1}G(s)^{C-1}&A_{-1}\hat{G}(s)
\end{array}\right].\end{equation} Observe that
\begin{equation}\small{\mathcal{D}\left[\begin{array}{cc} I& \hat{G}(s)^{C}
\\G(s)^{C}&I
\end{array}\right]=\left[\begin{array}{cc} A_{-1} & A_{-1} \hat{G}(s)^C\\0&0\\\vdots&\vdots \\0&0\\A_1 G(s)^C& A_1\end{array}\right]=-\mathcal{C}(s)\left[\begin{array}{cc} G(s) & \hat{G}(s)^{C-1} \\ G(s)^2 & \hat{G}(s)^{C-2} \\\vdots&\vdots\\G(s)^{C-2} & \hat{G}(s)^2\\ G(s)^{C-1} & \hat{G}(s)\end{array}\right]}
\end{equation}where we have used \eqref{eqG} and \eqref{eqGh}. This directly leads to $$(-\mathcal{C}(s))^{-1} \mathcal{D} \left[\begin{array}{cc} I& \hat{G}(s)^{C}
\\G(s)^{C}&I
\end{array}\right]=\left[\begin{array}{cc} G(s) & \hat{G}(s)^{C-1} \\ G(s)^2 & \hat{G}(s)^{C-2} \\\vdots&\vdots\\ G(s)^{C-1} & \hat{G}(s)\end{array}\right],$$which, pre-multiplied by the matrix $\mathcal{B}$, provides \eqref{toshow}. 
 \hfill\cqfd
\medskip

Observe that the expression (\ref{sol2}) for $\tilde{\vc R}_{k}(s)$ involves matrices of {at most} twice the size of the phase space {(such as $Z(s,C)$)}.
 The function $\vc R_k(t)$ is obtained by taking the inverse Laplace transform of $\tilde{\vc R}_{k}(s)$, which can be done numerically using the method of Abate and Whitt \cite{abate1995numerical} and the function $\mathsf{ilaplace.m}$ in Matlab.
 
 The next corollary provides us with the limit as $C\rightarrow\infty$ of the result stated in Proposition \ref{prop26}. It is a direct consequence of the fact that the matrices $G(s)$ and $ \hat{G}(s)$ are sub-stochastic for any $s>0$.
\begin{cor}\label{cor_prop26}
Assume that the series $\sum_{j=1}^{\infty} \hat{G}(s)^j H_0(s)\,\vc g_{k+j}(s)$ converges for any $k\geq 0$. Then, for any $s>0$, and $k\geq 0$, the limit as $C\rightarrow\infty$ of the Laplace transform of the expected reward function, conditional on the initial level $k$, is given by
\begin{equation}\label{sol2b}\tilde{\vc R}_k(s,\infty)=G(s)^k \vc v(s,\infty) +\vc\nu_k(s,\infty),\end{equation} where $$ \vc v(s,\infty)=-((B_0-sI)+A_1G(s))^{-1}(\vc g_0(s)+(B_0-sI)\,\vc\nu_0(s,\infty)+A_1\,\vc\nu_1(s,\infty)),$$and
\begin{eqnarray*}\vc\nu_k(s,\infty)&=&\sum_{j=0}^{k-1} G(s)^j H_0(s) \,\vc g_{k-j}(s)+\sum_{j=1}^{\infty} \hat{G}(s)^j H_0(s)\,\vc g_{k+j}(s),\quad k\geq 1.
\end{eqnarray*}\hfill\cqfd
\end{cor}
Note that we can write $\vc\nu_0(s,\infty)= \hat{G}(s) \vc\nu_1(s,\infty).$

\medskip
 
In addition to the blocks $\tilde{\vc R}_{k}(s)$ of the expected reward function, Proposition \ref{prop26} provides us with an explicit expression for the blocks $\tilde{D}_{k,\ell}(s)$ of the Laplace transform of the transient deviation matrix, as we show now.

\begin{prop}\label{proptdev} For $0\leq k\leq C$,
\begin{itemize}
\item if $\ell=0$,

$$\tilde{D}_{k,0}(s)=G(s)^k V(s,C) +\hat{G}(s)^{C-k}W(s,C)- \vc 1\vc\pi_0\,(1/s^2),$$  with
$$\left[\begin{array}{c} V(s,C)\\W(s,C)\end{array}\right]=(-s\,Z(s,C))^{-1}\left[\begin{array}{c} I\\0\end{array}\right],$$
\item if $1\leq \ell\leq C-1$, 
\begin{eqnarray*}\tilde{D}_{k,\ell}(s)&=&\left(G(s)^k V(s,C) +\hat{G}(s)^{C-k}W(s,C)\right.\\&&+\left.(1/s)G(s)^{k-\ell}\mathds{1}_{\{\ell\leq k\}}+(1/s)\hat{G}(s)^{\ell-k} \mathds{1}_{\{ \ell> k\}}\right)H_0(s)- \vc 1\vc\pi_\ell\,(1/s^2),\end{eqnarray*} with
$$\left[\begin{array}{c} V(s,C)\\W(s,C)\end{array}\right]=(-s\,Z(s,C))^{-1}\left[\begin{array}{c}(B_0-sI) \hat{G}(s)^\ell +A_1\hat{G}(s)^{\ell-1}\\A_{-1}G(s)^{C-1-\ell}+(C_0-sI)G(s)^{C-\ell}\end{array}\right],$$
\item if $\ell=C$,
$$\tilde{D}_{k,C}(s)=\left(G(s)^k V(s,C) +\hat{G}(s)^{C-k}(W(s,C)+(1/s)I)\right)H_0(s)- \vc 1\vc\pi_C\,(1/s^2),$$  with
$$\left[\begin{array}{c} V(s,C)\\W(s,C)\end{array}\right]=(-s\,Z(s,C))^{-1}\left[\begin{array}{c}(B_0-sI) \hat{G}(s)^C+A_1 \hat{G}(s)^{C-1} \\(C_0-A_0)-A_1 G(s)\end{array}\right],$$\end{itemize}
 where
$Z(s,C)$ is given by \eqref{ZsC}.

\end{prop} 

\noindent \textit{Proof.}
We use \eqref{LTTD} and Proposition \ref{prop26} with, for $0\leq \ell\leq C$ and $1\leq i\leq n$, $\vc g_\ell=\vc e_i$ and $\vc g_j=\vc 0$ for $j\neq \ell$. 
\hfill\cqfd

\medskip

The limit as $C\rightarrow\infty$ of Proposition \ref{proptdev} again follows easily, {and provides us with an analytical expression for the blocks of the infinite matrix $\tilde{D}(s)$ corresponding to a QBD with no upper bound on the levels.}
\begin{cor} For $k\geq 0$,
\begin{itemize}
\item if $\ell=0$,
$$\tilde{D}_{k,0}(s,\infty)=G(s)^k (-s)^{-1}((B_0-sI)+A_1G(s))^{-1} - \vc 1\vc\pi_0\,(1/s^2),$$ 
\item if $\ell\geq 1$, 
\begin{eqnarray*}\tilde{D}_{k,\ell}(s,\infty)&=&\left(G(s)^k V(s,\infty)+(1/s)G(s)^{k-\ell}\mathds{1}_{\{\ell\leq k\}}\right.\\&&\left. +(1/s)\hat{G}(s)^{\ell-k} \mathds{1}_{\{ \ell> k\}}\right)H_0(s)- \vc 1\vc\pi_\ell\,(1/s^2),\end{eqnarray*} with
$$V(s,\infty)=(-s)^{-1}((B_0-sI)+A_1G(s))^{-1}((B_0-sI) \hat{G}(s) +A_1)\hat{G}(s)^{\ell-1}.$$\hfill\cqfd
\end{itemize}
\end{cor}

 \subsection{Asymptotic regime}

Next, we concentrate on the asymptotic properties of the expected reward function ${\vc R}_k(t)$ for large values of $t$. By decomposing \eqref{lin_ass} into blocks, we have
$$\bar{\vc R}_k(t)=\sum_{0\leq \ell\leq C}((\vc\pi_\ell\,\vc g_\ell)\vc 1\,t\, +D_{k,\ell}\,\vc g_\ell).$$
We determine now an explicit expression for the blocks $D_{k,\ell}$ of the deviation matrix for $0\leq k,\ell\leq C$,
   together with the mean first passage times to any level $\ell$ in the QBD process. Our method is based on the relationship between the entries of the deviation matrix and mean first passage times, and involves the solution of finite systems of matrix difference equations similar to the ones we solved in Section \ref{tr}. Indeed, we can write
$$D_{(k,i)(\ell,j)}=\pi_{(\ell,j)}\,[M_{\vc\pi(\ell,j)}-M_{(k,i)(\ell,j)}],$$ where $M_{(k,i)(\ell,j)}$ is the mean first entrance time to $(\ell,j)$ from $(k,i)$, and $M_{\vc\pi(\ell,j)}$ is the mean first entrance time to $(\ell,j)$ if the state at time 0 has the stationary distribution $\vc\pi$. Define the block-matrices $M_{k,\ell}=\left(M_{(k,i)(\ell,j)}\right)_{1\leq i,j\leq n}$, for $0\leq k,\ell\leq C$. In matrix form, we have 
\begin{equation}\label{dev_formula}D_{k,\ell}=\left[\left(\vc 1_n\otimes \sum_{0\leq x\leq C} \vc \pi_x \,M_{x, \ell}\right)-M_{k,\ell}\right]\mbox{diag}(\vc\pi_\ell).\end{equation}

Let us fix level $\ell$ and define $\vc m^{(j)}_{k,\ell}=  M_{k,\ell}\,\vc e_j$, the $j$th column of $M_{k,\ell}$, for $1\leq j\leq n$ and $0\leq k\leq C$. So  $(\vc m^{(j)}_{k,\ell})_i=M_{(k,i)(\ell,j)}$.  Since the diagonal of $M_{\ell,\ell}$ is null we must have $(\vc m^{(j)}_{\ell,\ell})_j=0$. We introduce the notation $\bar{A}_{-1}^{(j)}$ and $\bar{A}_{1}^{(j)}$ for the matrices obtained by replacing the $j$th row in $A_{-1}$ and $A_{1}$, respectively, by $\vc 0^\top$, and $\bar{B}_{0}^{(j)} $, $\bar{A}_{0}^{(j)} $, $\bar{C}_{0}^{(j)} $ for the matrices obtained by replacing the $j$th row in $B_0$, $A_0$ and $C_0$, respectively, by $-\vc e_j^\top$. 

\begin{prop}\label{mj} For any fixed level $0\leq \ell\leq C$ and any phase $1\leq j\leq n$, the vectors $\vc m^{(j)}_{k,\ell}$ satisfy the system of {matrix second-order difference equations}
\begin{eqnarray} \label{bccc}A_{-1}\,\vc m^{(j)}_{k-1,\ell}+A_{0}\,\vc m^{(j)}_{k,\ell}+A_1 \,\vc m^{(j)}_{k+1,\ell}
&=&-\vc 1\end{eqnarray} for $1\leq k \leq \ell-1$ and  $\ell+1\leq k\leq C-1$, with the following boundary conditions, depending on the value of $\ell$:
\begin{itemize}\item for $\ell=0$,
\begin{eqnarray}\label{bc2b}\bar{B}_0^{(j)}\,\vc m^{(j)}_{0,\ell}+\bar{A}_1^{(j)}\,\vc m^{(j)}_{1,\ell}&=&-\vc 1+\vc e_j,\\\label{bc1b}A_{-1}\,\vc m^{(j)}_{C-1,\ell}+C_{0}\,\vc m^{(j)}_{C,\ell}
&{=}&{-\vc 1},\end{eqnarray} \item for $1\leq \ell\leq C-1$,
\begin{eqnarray}\label{bc1}B_0\,\vc m^{(j)}_{0,\ell}+A_1\,\vc m^{(j)}_{1,\ell}&=&-\vc 1,\\\label{bc2}
\bar{A}_{-1}^{(j)}\vc m^{(j)}_{\ell-1,\ell}+\bar{A_0}^{(j)} \vc m^{(j)}_{\ell,\ell}+\bar{A}_{1}^{(j)}\vc m^{(j)}_{\ell+1,\ell}&=&-\vc 1+\vc e_j\\\label{last}A_{-1}\,\vc m^{(j)}_{C-1,\ell}+C_{0}\,\vc m^{(j)}_{C,\ell}
&{=}&{-\vc 1},\end{eqnarray}
 \item for $\ell=C$,
\begin{eqnarray}\label{bc1c}B_0\,\vc m^{(j)}_{0,\ell}+A_1\,\vc m^{(j)}_{1,\ell}&=&-\vc 1,\\\label{bc2c}\bar{A}_{-1}^{(j)}\,\vc m^{(j)}_{C-1,\ell}+\bar{C}_{0}^{(j)}\,\vc m^{(j)}_{C,\ell}
&{=}&{-\vc 1+\vc e_j}.\end{eqnarray}\end{itemize}

\end{prop} 

\noindent \textit{Proof.} For a fixed value of $\ell$, let $m^{(j)}_{(k,i)}=(\vc m^{(j)}_{k,\ell})_i$, for $1\leq i\leq n$. By conditioning on the epoch where the process first leaves the state $(k,i)$, we obtain for $1\leq k\leq C-1$, $\ell+1\leq k\leq C-1$, and $1\leq j\leq n$, 
\begin{eqnarray}\nonumber m^{(j)}_{(k,i)}&=&\dfrac{1}{(-A_0)_{ii}}+\sum_{x\neq i}\dfrac{(A_0)_{ix}}{(-A_0)_{ii}} m^{(j)}_{(k,\ell)}
\\\label{gen1}&&+\sum_{x}\dfrac{(A_1)_{ix}}{(-A_0)_{ii}} m^{(j)}_{(k+1,x)}+\sum_{x}\dfrac{(A_{-1})_{ix}}{(-A_0)_{ii}} m^{(j)}_{(k-1,x)},\end{eqnarray}which in matrix form gives \eqref{bccc}. 

A similar argument leads to the boundary equations \eqref{bc1b}, \eqref{bc1}, \eqref{last}, and \eqref{bc1c}. 
Finally, the boundary equations \eqref{bc2b}, \eqref{bc2}, and \eqref{bc2c} are obtained by adding the constraint that when $k=\ell$, $m^{(j)}_{(\ell,j)}=0$. \hfill\cqfd

\medskip


For $0\leq k\leq C$, define the vectors
\begin{equation}\label{sigk}\vc\mu_k(C)=\sum_{j=0}^{k-1} G^j H_0 \,\vc 1+\sum_{j=1}^{C-k} \hat{G}^j H_0\,\vc 1. \end{equation} The next theorem provides an explicit expression for the columns of the mean first passage time matrices $M_{k,\ell}$. 

\begin{prop}\label{prop2} In the non null-recurrent case, for any $0\leq k\leq C$, the vector $\vc m^{(j)}_{k,\ell}$ has the following explicit expression, depending on the value of $\ell$:
\begin{itemize}\item if $\ell=0$,
$$\vc m^{(j)}_{k,0}=G^k \vc v^{(j)}+\hat{G}^{C-k} \vc w^{(j)} + \vc \mu_k(C),$$ where $\vc \mu_k(C)$ is given by \eqref{sigk} and
$$\left[\begin{array}{c} \vc v^{(j)}\\ \vc w^{(j)}\end{array}\right]=(-Z^{(j)}(C))^{-1}\left[\begin{array}{c}\vc 1-\vc e_j+\bar{B}_0^{(j)} \vc \mu_0(C) + \bar{A}_1^{(j)}\vc \mu_1(C)\\\vc 1+{A}_{-1}\vc\mu_{C-1}(C)+{C}_{0}\vc \mu_C(C)\end{array}\right]$$with
\begin{eqnarray*}Z^{(j)}(C)&=&\left[\begin{array}{cc} \bar{B}_0^{(j)}+\bar{A}_1^{(j)}G& (\bar{B}_0^{(j)}\hat{G}+\bar{A}_1^{(j)})\hat{G}^{C-1}
\\({A}_{-1}+{C}_{0}G) G^{C-1}& {A}_{-1}\hat{G}+{C}_{0}
\end{array}\right]\\&=&\mathring{Q}^{(j)}(C) \left[\begin{array}{cc} I& \hat{G}^{C}
\\G^{C}&I
\end{array}\right],\end{eqnarray*}where $\mathring{Q}^{(j)}(C)$ is the generator of the transient Markov chain obtained from the QBD process in which absorption can happen from state $(0,j)$ at rate $1$, and restricted to levels 0 and $C$.
\item if $1\leq \ell\leq C-2$,
\begin{eqnarray*}\vc m^{(j)}_{k,\ell}&=&\left(G^k \vc v^{-(j)}+\hat{G}^{\ell-k} \vc w^{-(j)}\right) \mathds{1}_{\{k\leq \ell\}}\\&&+\left(G^{k-\ell-1} \vc v^{+(j)}+\hat{G}^{C-k} \vc w^{+(j)}\right)\mathds{1}_{\{k\geq  \ell+1\}} \\&& + \vc \mu _k(C),\end{eqnarray*}
where $\vc \mu_k(C)$ is given by \eqref{sigk} and
$$\small{\left[\begin{array}{c} \vc v^{-(j)}\\ \vc w^{-(j)}\\\vc v^{+(j)}\\\vc w^{+(j)}\end{array}\right]=(-Z^{(j)}(\ell,C))^{-1}\left[\begin{array}{c}\vc 1+B_0 \vc \mu_0(C) + A_1\vc \mu_1(C)\\\vc 1-\vc e_j+\bar{A}_{-1}^{(j)}\vc\mu_{\ell-1}(C)+\bar{A}_{0}^{(j)}\vc \mu_\ell(C)+\bar{A}_{1}^{(j)}\vc\mu_{\ell+1}(C)\\\vc 1+A_{-1}\vc\mu_{\ell}(C)+A_0\vc\mu_{\ell+1}(C)+A_{1}\vc\mu_{\ell+2}(C)\\ \vc 1+ A_{-1}\vc\mu_{C-1}(C)+C_0\vc\mu_C(C)\end{array}\right]}$$with
\begin{eqnarray*}\lefteqn{Z^{(j)}(\ell,C)}\\&=&\small{\left[\begin{array}{cccc} B_0+A_1G& (B_0\hat{G}+A_1)\hat{G}^{\ell-1}& 0 & 0\\(\bar{A}_{-1}^{(j)}+\bar{A}_{0}^{(j)}G) G^{\ell-1}&\bar{A}_{-1}^{(j)}\hat{G}+\bar{A}_{0}^{(j)}& \bar{A}_1^{(j)} &\bar{A}_1^{(j)} \hat{G}^{C-\ell-1}\\A_{-1} G^\ell & A_{-1} & A_0+A_1G & (A_0 \hat{G}+A_1)\hat{G}^{C-\ell-2}\\0 &0&(C_0G+A_{-1})G^{C-\ell-2} & C_0+A_{-1} \hat{G}
\end{array}\right]}\\\normalsize &=&\mathring{Q}^{(j)}(\ell,C)\left[\begin{array}{cccc} I& \hat{G}^{\ell}& 0 &0
\\G^{\ell}&I & 0 &0\\0 & 0 & I & \hat{G}^{C-\ell-1}\\ 0&0&G^{C-\ell-1} & I
\end{array}\right],\end{eqnarray*}
where $\mathring{Q}^{(j)}(\ell,C)$ is the generator of the 
 transient Markov chain obtained from the QBD process  in which absorption can happen from state $(\ell,j)$ at rate $1$, and restricted to levels 0, $\ell$, $\ell+1$, and $C$.
\item if $\ell=C-1$,
$$\vc m^{(j)}_{k,C-1}=\left(G^k \vc v^{(j)}+\hat{G}^{C-1-k} \vc w^{(j)}\right)\mathds{1}_{\{k\leq C-1\}}+\vc x^{(j)}\mathds{1}_{\{k= C\}} + \vc \mu_k(C),$$ where $\vc \mu_k(C)$ is given by \eqref{sigk} and
$$\small{\left[\begin{array}{c} \vc v^{(j)}\\ \vc w^{(j)}\\\vc x^{(j)}\end{array}\right]=(-Z^{(j)}(C))^{-1}\left[\begin{array}{c}\vc 1+B_0 \vc \mu_0 (C)+ A_1\vc \mu_1(C)\\\vc 1-\vc e_j+\bar{A}_{-1}^{(j)}\vc\mu_{C-2}(C)+\bar{A}_{0}^{(j)}\vc \mu_{C-1}(C)+\bar{A}_{1}^{(j)}\vc \mu_{C}(C)\\\vc 1+A_{-1}\vc \mu_{C-1}(C)+C_0\vc \mu_{C}(C)\end{array}\right]}$$with
\begin{eqnarray*}Z^{(j)}(C)&=&\left[\begin{array}{ccc} B_0+A_1G& (B_0\hat{G}+A_1)\hat{G}^{C-2}& 0
\\(\bar{A}_{-1}^{(j)}+\bar{A}_{0}^{(j)}G) G^{C-2}& \bar{A}_{-1}^{(j)}\hat{G}+\bar{A}_{0}^{(j)} & \bar{A}_{1}^{(j)}\\ A_{-1} G^{C-1} & A_{-1}& C_0
\end{array}\right]\\&=&\mathring{Q}^{(j)}(C) \left[\begin{array}{ccc} I& \hat{G}^{C-1}& 0
\\G^{C-1}&I & 0\\0&0 & I
\end{array}\right],\end{eqnarray*}where $\mathring{Q}^{(j)}(C)$ is the generator of the transient Markov chain obtained from the QBD process in which absorption can happen from state $(C-1,j)$ at rate $1$, and restricted to levels 0, $C-1$ and $C$.
\item if $\ell=C$,
$$\vc m^{(j)}_{k,C}=G^k \vc v^{(j)}+\hat{G}^{C-k} \vc w^{(j)} + \vc \mu_k(C),$$ where $\vc \mu_k(C)$ is given by \eqref{sigk} and
$$\left[\begin{array}{c} \vc v^{(j)}\\ \vc w^{(j)}\end{array}\right]=(-Z^{(j)}(C))^{-1}\left[\begin{array}{c}\vc 1+B_0 \vc \mu_0(C) + A_1\vc \mu_1(C)\\\vc 1-\vc e_j+\bar{A}_{-1}^{(j)}\vc\mu_{C-1}(C)+\bar{C}_{0}^{(j)}\vc \mu_C(C)\end{array}\right]$$with
\begin{eqnarray*}Z^{(j)}(C)&=&\left[\begin{array}{cc} B_0+A_1G& (B_0\hat{G}+A_1)\hat{G}^{C-1}
\\(\bar{A}_{-1}^{(j)}+\bar{C}_{0}^{(j)}G) G^{C-1}& \bar{A}_{-1}^{(j)}\hat{G}+\bar{C}_{0}^{(j)}
\end{array}\right]\\&=&\mathring{Q}^{(j)}(C) \left[\begin{array}{cc} I& \hat{G}^{C}
\\G^{C}&I
\end{array}\right],\end{eqnarray*}where $\mathring{Q}^{(j)}(C)$ is the generator of the transient Markov chain obtained from the QBD process in which absorption can happen from state $(C,j)$ at rate $1$, and restricted to levels 0 and $C$.
\end{itemize}
\end{prop} 

\noindent\textit{ Proof.} The result follows from Lemma \ref{solhomu} applied to the system of difference equations \eqref{bccc}. The arbitrary vectors are then determined using the boundary conditions given in Proposition \ref{mj}.

 When $\ell=0$ and $\ell=C$, there are two boundary equations, determining two arbitrary vectors. 
 
 When $1\leq \ell \leq C-2$, the solution depends on whether $0\leq k\leq \ell$ or $\ell+1\leq k\leq C$. There are four boundary conditions, namely the three equations described in Proposition \ref{mj}, in addition to one boundary equation {obtained by taking $k=\ell+1$ in \eqref{bccc}}. This determines the four arbitrary vectors. 
 
 When $\ell=C-1$ the solution depends on whether $0\leq k\leq C-1$ or $k=C$. There are three boundary equations, as described in Proposition \ref{mj}, which determine the three arbitrary vectors. 

The decomposition of the matrices $Z^{(j)}(C)$ and $Z^{(j)}(\ell,C)$ into the product of a non-conservative generator and the matrices involving $G$ and $\hat{G}$ follows from the same arguments as those used in the proof of Proposition \ref{prop26}. {Absorption from state $(\ell,j)$ at rate 1 comes from the definition of $\bar{B}_{0}^{(j)} $, $\bar{A}_{0}^{(j)} $, and $\bar{C}_{0}^{(j)}$}. The matrices involving $G$ and $\hat{G}$ are invertible because
sp$(G)<1$ or sp$(\hat{G})<1$ in the non null-recurrent case. \hfill\cqfd

\medskip

\medskip

Observe that in the transient case the series $\sum_{j=1}^{\infty} \hat{G}^j H_0\,\vc 1$ diverges, and the limit as $C\rightarrow\infty$ of $\vc\mu_k(C)$ given in \eqref{sigk} is infinite, while in the positive recurrent case, \begin{equation}\label{sigklim}\lim_{C\rightarrow\infty}\vc\mu_k(C)=\vc\mu_k(\infty)=\sum_{j=0}^{k-1} G^j H_0 \,\vc 1+((I-\hat{G})^{-1}-I) H_0\,\vc 1.\end{equation} This leads to the following corollary to Proposition \ref{prop2} for the mean first passage times in an unrestricted positive recurrent QBD process:
\begin{cor} In the positive recurrent case, for any $k\geq 0$, depending on the value of $\ell$, the vector $\vc m^{(j)}_{k,\ell}$ has the explicit expression:
\begin{itemize}\item if $\ell=0$,
$$\vc m^{(j)}_{k,0}(\infty)=G^k \vc v^{(j)}(\infty)+ \vc \mu_k(\infty),$$ where $\vc\mu_k(\infty)$ is given by \eqref{sigklim}, and
$$\vc v^{(j)}(\infty)=-( \bar{B}_0^{(j)}+\bar{A}_1^{(j)}G)^{-1}(\vc 1-\vc e_j+\bar{B}_0^{(j)} \vc \mu_0(\infty) + \bar{A}_1^{(j)}\vc \mu_1(\infty))$$
\item if $\ell\geq 1$,
$$\vc m^{(j)}_{k,\ell}=\left(G^k \vc v^{-(j)}+\hat{G}^{\ell-k} \vc w^{-(j)} \right)\mathds{1}_{\{k\leq \ell\}}+G^{k-\ell-1} \vc v^{+(j)}\mathds{1}_{\{k\geq  \ell+1\}}  + \vc \mu _k(\infty),$$
where $\vc\mu_k(\infty)$ is given by \eqref{sigklim} and
$$\small{\left[\begin{array}{c} \vc v^{-(j)}\\ \vc w^{-(j)}\\\vc v^{+(j)}\end{array}\right]=(-W^{(j)}(\ell))^{-1}\left[\begin{array}{c}\vc 1+B_0 \vc \mu_0(\infty) + A_1\vc \mu_1(\infty)\\\vc 1-\vc e_j+\bar{A}_{-1}^{(j)}\vc\mu_{\ell-1}(\infty)+\bar{A}_{0}^{(j)}\vc \mu_\ell(\infty)+\bar{A}_{1}^{(j)}\vc\mu_{\ell+1}(\infty)\\\vc 1+A_{-1}\vc\mu_{\ell}(\infty)+A_0\vc\mu_{\ell+1}(\infty)+A_{1}\vc\mu_{\ell+2}(\infty)\end{array}\right]}$$with
\begin{eqnarray*}W^{(j)}(\ell)&=&\left[\begin{array}{ccc} B_0+A_1G& (B_0\hat{G}+A_1)\hat{G}^{\ell-1}& 0 \\(\bar{A}_{-1}^{(j)}+\bar{A}_{0}^{(j)}G) G^{\ell-1}&\bar{A}_{-1}^{(j)}\hat{G}+\bar{A}_{0}^{(j)}& \bar{A}_1^{(j)}\\A_{-1} G^\ell & A_{-1} & A_0+A_1G 
\end{array}\right].\end{eqnarray*}
\end{itemize}\hfill\cqfd
\end{cor}

\medskip
To complete the characterisation of the block matrices $D_{k,\ell}$ using \eqref{dev_formula}, it remains for us to compute the blocks $\vc\pi_x$ of the stationary distribution of the QBD process for $0\leq x\leq C$. This can be done following Theorem 10.3.2 in \cite{latram} or Hajek \cite{hajek}, adapted to the continuous-time setting. This involves the matrices $R$ and $\hat{R}$ of rates of sojourn in level $k+1$, respectively $k-1$, per unit of the local time in level $k$ in the corresponding unrestricted QBD process. These matrices can be expressed in terms of $G$ and $\hat{G}$ as $R=A_1(-(A_0+A_1G))^{-1}$ and $\hat{R}=A_{-1}(-(A_0+A_{-1}\hat{G}))^{-1}$, see \cite{latram}.

\begin{cor} In the non null-recurrent case, the stationary distribution of the QBD process with  transition matrix \eqref{GeneratorQ} is given by
$$\vc\pi_k=\vc v_0\, R^k+\vc v_C\, \hat{R}^{C-k},\quad 0\leq k\leq C,$$ where $(\vc v_0,\vc v_C)$ is the solution of the system
$$(\vc v_0,\vc v_C)\,\left[\begin{array}{cc}B_0+R\,A_{-1}& R^{C-1}\,(R C_0+A_1)\\\hat{R}^{C-1}\,(\hat{R} B_0+A_{-1}) & C_0+\hat{R}\,A_1
\end{array}\right]=\vc 0,$$ and $$ \vc v_0 \sum_{0\leq i\leq C} R^i\,\vc 1+\vc v_C \sum_{0\leq i\leq C} \hat{R}^i\,\vc 1=1.$$
\end{cor}

\section{A perturbation theory approach}\label{pta}

In this section, we answer Question \emph{(ii)} on Page~\pageref{questiondeux} and we provide an
expression for the deviation matrix of a QBD with maximal level $C$ in terms of the deviation matrix of the same QBD with capacity $C-1$, in transient and asymptotic regimes. To highlight the number of levels in the QBDs, we shall denote by $Q^{(C)}$ the generator \eqref{GeneratorQ} of a QBD with maximum level $C$, and by $D^{(C)}(t)$, {$\tilde{D}^{(C)}(s)$,} and $D^{(C)}$ the corresponding deviation matrices, for $C\geq 1$.

The key idea behind the derivation of the recursions is the decomposition of the generator $Q^{(C)}$ into a sum of another generator and a perturbation matrix, so that $Q^{(C)}$ can  be seen
as a block-element updating problem, the block version of the element-updating
problem described in Langville and Meyer \cite{langville2006updating}.  Precisely, we write
\begin{equation}
Q^{(C)}={T}^{(C)}+E_{C-1}^{(C)}\Delta^{(C)},\label{SumDecompo}
\end{equation}
where
\begin{eqnarray}\label{eq:Qtilde}
{T}^{(C)}&=&\left[\begin{array}{cccc|c}
 &  &  &  & 0\\
 &  & Q^{(C-1)} &  & \vdots\\
 &  &  &  & 0\\
\hline 0 & \cdots & 0 & A_{-1} & C_{0}
\end{array}\right],\\E_{C-1}^{(C)} & =&\left[\begin{array}{ccccc}
0 & \cdots & 0 & I & 0\end{array}\right]^{\top},\nonumber \\\label{eq:delta}
\Delta^{(C)}&= & [\begin{array}{cccccc}
0 & \cdots & 0 &  & A_{0}-C_{0} & A_{1}\end{array}],
\end{eqnarray}where $E_{C-1}^{(C)} $ has dimensions $n(C+1) \times n$ and $\Delta^{(C)}$ has dimensions $n \times n(C+1)$.
Observe that ${T}^{(C)}$ is the generator of a reducible transient structured Markov chain, and its stationary distribution is given by
$$\boldsymbol{\phi}=[\;\boldsymbol{\pi}^{(C-1)},\;\vc 0\;],$$ where
$\boldsymbol{\pi}^{(C-1)}$ is the stationary distribution of $Q^{(C-1)}$. 
The group inverse of $T^{(C)}$
exists by \cite[Theorem 8.5.5]{campbell2009generalized},
 and may
be expressed as a function of $\boldsymbol{\pi}^{(C-1)}$ and 
the deviation matrix $D^{(C-1)}=-(Q^{(C-1)})^{\#} $, as shown in the next lemma.

\begin{lem}\label{devTC}The group inverse of $T^{(C)}$ is given by
\begin{equation}
(T^{(C)})^{\#}  =\left[\begin{array}{c|c}
{-D^{(C-1)}}  & 0\\\hline
C_{0}^{-1}\,M^{(C)}\,D^{(C-1)}+M^{(C)}\,\boldsymbol{1\pi}^{(C-1)} & C_{0}^{-1}
\end{array}\right],\label{eq:QsharpeTilde}\end{equation}where $M^{(C)}$ is an $n \times n C$  matrix defined as \begin{equation}
M^{(C)} =\left[\begin{array}{cccc}
0 & \cdots & 0 & A_{-1}\end{array}\right].  \label{VectM}
\end{equation}
\end{lem}

\noindent \textit{Proof.} We check that $(T^{(C)})^{\#}$
is the group inverse of $T^{(C)}$ by direct verification of (\ref{Cesaro}, \ref{Cesaro2}).\hfill{}\cqfd


\medskip
 
First, we use the decomposition \eqref{SumDecompo} to derive an expression for $\boldsymbol{\pi}^{(C)}$ in terms of $\boldsymbol{\pi}^{(C-1)}$. 

\begin{lem}\label{pirec} The stationary distribution of $Q^{(C)}$ may be expressed in terms of the stationary distribution of $Q^{(C-1)}$ as
\[
\boldsymbol{\pi}^{(C)}=[\;\boldsymbol{\pi}^{(C-1)},\;\vc 0\;](I+E_{C-1}^{(C)}\Delta^{(C)} (T^{(C)})^{\#})^{-1}.
\]
\end{lem}

\noindent \textit{Proof.}  Our argument is a slight modification of the proof of Rising \cite[Lemma 5.1]{rising1991applications}, which deals with irreducible Markov matrices.  Let $A$ and $\widetilde A$ be two generators on the same state space, assume that each has a single irreducible class of states, not necessarily the same, so that each has a unique stationary probability vector.  Denote the stationary distribution vectors as $\vc\alpha$ and $\widetilde{\vc\alpha}$, respectively.

By \eqref{eq:deveq2}, $\widetilde{\vc\alpha} A A^\# = \widetilde{\vc\alpha} - \vc\alpha$ and so
\[
\widetilde{\vc\alpha} (I + (\widetilde A - A) A^\#) = \vc\alpha
\]
since $\widetilde{\vc\alpha} \widetilde A = \vc 0$.  If $A$ and $\widetilde A$ are irreducible, it results from \cite[Lemma 5.1]{rising1991applications} that $I + (\widetilde A - A) A^\#$ is non-singular.  A key argument in \cite{rising1991applications} is that the kernel of both $A$ and $\widetilde A$ is limited to Span($\{\vc 1\}$), and so we repeat it verbatim to conclude that $I + (\widetilde A - A) A^\#$  is non-singular in our case as well.
Thus,
\begin{equation}
\widetilde{\vc\alpha} = \vc\alpha (I + (\widetilde A - A) A^\#)^{-1}
   \label{StationPerturbed}
\end{equation}
and the lemma is proved, once we set  $A=T^{(C)}$ and $\widetilde A = Q^{(C)}$.  
\hfill{}\cqfd
\medskip

We thus see that $\boldsymbol{\pi}^{(C)}$ can be expressed in terms of $\boldsymbol{\pi}^{(C-1)}$ and $D^{(C-1)}$ through $(T^{(C)})^{\#}$.
In the next two sections, we use the decomposition \eqref{SumDecompo} to derive  recursive formulae for {$\tilde{D}^{(C)}(s)$} and $D^{(C)}$.

\subsection{Time-dependent deviation matrix}
Recall that the Laplace transform of the transient deviation matrix of $Q^{(C)}$ is given by
\begin{equation}\label{trans}\tilde{D}^{(C)}(s)= (1/s) \left(sI-Q^{(C)}\right)^{-1}-(1/s^2)\vc 1\vc\pi^{(C)}. \end{equation} We can express $\tilde{D}^{(C)}(s)$ in terms of $\tilde{D}^{(C-1)}(s)$ thanks to the fact that $(sI-Q^{(C)})^{-1}$ can be expressed in terms
of $(sI-Q^{(C-1)})^{-1}$ via the Sherman-Morrison-Woodbury (SMW) formula, which we restate here.

Let $M$ be a finite non-singular matrix of order $m$ and let $U$ and $V$ be matrices of dimension $m \times k$ and $k \times m$, respectively.   If $M+UV^\top$ is
non-singular, then $I+V^\top M^{-1}U$ is non-singular as well  and  
\begin{equation}
(M+UV^\top)^{-1}=M^{-1}\left(I-U(I+V^\top M^{-1}U)^{-1}V^\top M^{-1}\right).\label{eq:SMW}
\end{equation}
As a direct consequence of this formula, we have the following

\begin{cor} The inverse of $(sI-Q^{(C)})$ satisfies the recursive equation
\begin{eqnarray*}\lefteqn{(sI-Q^{(C)})^{-1}}\\&&=(sI-T^{(C)})^{-1} \left(I+E_{C-1}^{(C)}(I-\Delta^{(C)} (sI-T^{(C)})^{-1} E_{C-1}^{(C)})^{-1}\Delta^{(C)} (sI-T^{(C)})^{-1}\right),\end{eqnarray*}where
$\Delta^{(C)}$ is given in \eqref{eq:delta}, and
\begin{eqnarray*}(sI-T^{(C)})^{-1}&=&\left[\begin{array}{c|c} I & 0\\\hline (sI-C_0)^{-1}M^{(C)}& (sI-C_0)^{-1}
\end{array}\right]\left[\begin{array}{c|c}(sI-Q^{(C-1)})^{-1}&0\\\hline0& I\end{array}\right],\end{eqnarray*} where $M^{(C)}$ is given in \eqref{VectM}.
\end{cor} 

\noindent \textit{Proof.}  We use the decomposition \eqref{SumDecompo} and apply \eqref{eq:SMW} with $M=sI-T^{(C)}$, $U=E_{C-1}^{(C)}$, and $V=-\Delta^{(C)\top}$. The expression for $(sI-T^{(C)})^{-1}$ can be checked by verifying that $(sI-T^{(C)})^{-1}(sI-T^{(C)})=I$ using \eqref{eq:Qtilde}. 
%
%
\hfill{}\cqfd

\medskip
The inverse $(sI-Q^{(C)})^{-1}$ is expressed in terms
of $(sI-Q^{(C-1)})^{-1}$ and, by \eqref{trans}, 
\[
(sI-Q^{(C-1)})^{-1}=s\tilde{D}^{(C-1)}(s)+(1/s)\vc 1\vc\pi^{(C-1)}.
\]
Furthermore, since $\boldsymbol{\pi}^{(C)}$ can be expressed in terms of $\boldsymbol{\pi}^{(C-1)}$ and $D^{(C-1)}$ by Lemmas \ref{devTC} and \ref{pirec}, we conclude that $\tilde{D}^{(C)}(s)$ may be expressed in terms of $\tilde{D}^{(C-1)}(s)$, $\boldsymbol{\pi}^{(C-1)}$, and $D^{(C-1)}$. The recursive computation of the transient deviation matrices $\tilde{D}^{(C)}(s)$ may therefore be done together with the recursive computation of the stationary distribution vectors $\boldsymbol{\pi}^{(C)}$ and the deviation matrices ${D}^{(C)}$, for $C\geq 1$.

\subsection{Deviation matrix}

Now, we focus on obtaining an expression for $D^{(C)}$ in terms of $D^{(C-1)}$ using the decomposition \eqref{SumDecompo}.
We start from Theorem 3.1 in
  \cite{langville2006updating},
 which gives
 an expression
   for the group inverse of a
   one-element updated Markov chain generator, and we extend it to a one-block update. Starting from a generator $Q$, we consider
a new generator $\tilde{Q}=Q+A$ where $A$ is
a matrix with one non-zero block-row only.   We denote
by $E_{K}^{(C)}=[\begin{array}{cccccc}
0 & \cdots & 0 & I & \cdots & 0\end{array}]^\top$ the block-column with identity as $K$-th block and $C$ zero blocks
elsewhere, and by $P$ a block-row vector containing $C+1$ blocks, then $A = E_{K}^{(C)}P$ for some $K$. 

\begin{prop}\label{QSharpeExpress} Let $Q$ be the generator of
a Markov process with deviation matrix $D$ and with stationary distribution
$\boldsymbol{q}$. Suppose that $\tilde{Q}=Q+E_{K}^{(C)}P$  is
the generator of an irreducible Markov process. The deviation matrix
$\tilde{D}$ of $\tilde{Q}$ is given by 
\begin{eqnarray}
\tilde{D}&=&\left(I-\boldsymbol{1\tilde{\pi}}\right)D\left(I-E_{K}^{(C)}PD\right)^{-1}\label{QsharpeExpr}\\&=&\left(I-\boldsymbol{1\tilde{\pi}}\right)D \left(I+E_{K}^{(C)}(I-PDE_{K}^{(C)})^{-1}PD\right), \label{QsharpeExpr2}
\end{eqnarray}
where $\tilde{\boldsymbol{\pi}}$ is the stationary distribution of
$\tilde{Q}$. \end{prop}

\noindent \textit{Proof.} Observe first that the inverse in \eqref{QsharpeExpr} is well defined by \cite[Lemma 5.1]{rising1991applications}. We first show \eqref{QsharpeExpr} by direct verification of (\ref{eq:deveq2}, \ref{eq:deveq1}). Replacing $\tilde{D}$ with the right-hand side of  (\ref{QsharpeExpr}), we obtain 
\begin{align*}
\tilde{Q}\tilde{D} & =\tilde{Q}D\left(I-E_{K}^{(C)}PD\right)^{-1}  \\
\intertext{\mbox{since }\ensuremath{\tilde{Q}\boldsymbol{1}=\vc 0,}} 
& =(\vc 1\boldsymbol{q}-I)\left(I-E_{K}^{(C)}PD\right)^{-1}+E_{K}^{(C)}PD\left(I-E_{K}^{(C)}PD\right)^{-1}
\\
\intertext{\mbox{as }\ensuremath{\tilde{Q}={Q}+E_{K}^{(C)}P}\mbox{ and }
\ensuremath{{Q}D=\vc 1 \boldsymbol{q}-I}
by (\ref{eq:deveq2})\mbox{, } }
&
 =\vc 1\tilde{\boldsymbol{\pi}}-I,
\end{align*}
since $\boldsymbol{q}(I-E_{K}^{(C)}PD)^{-1}=\tilde{\boldsymbol{\pi}},$
by (\ref{StationPerturbed}), where $A$ is replaced by $Q$ and $\widetilde A$ by $Q+E_{K}^{(C)}P$.
Next, 
\begin{equation}
    \label{Simpl3}
\left(I-E_{K}^{(C)}PD\right)^{-1} 
  =I+E_{K}^{(C)}(I-PDE_{K}^{(C)})^{-1}PD
\end{equation}
by the SMW formula (\ref{eq:SMW}), and injecting (\ref{Simpl3}) into (\ref{QsharpeExpr}) 
provides (\ref{QsharpeExpr2}).
\hfill{}\cqfd

\medskip

\begin{rem}  \rm
Equation  (\ref{QsharpeExpr2}) can be seen as an extension of the SMW formula (\ref{eq:SMW}) to
singular matrices.
\end{rem}

\begin{rem}  \rm
As observed in \cite{langville2006updating}, the cost of straightforward computation of (\ref{QsharpeExpr}) may be high
if one updates more than one row, as it requires the inversion of the  matrix $I-E_{K}^{(C)}PD$ whose size is the same as $\tilde{Q}$. Equation (\ref{QsharpeExpr2}) reduces the computational cost because the size of $I-PDE_{K}^{(C)}$ is just one block.
\end{rem}

Since the decomposition of $Q^{(C)}$ given in (\ref{SumDecompo})
is a block-element-updating of the matrix $T^{(C)}$ by 
$E_{C-1}^{(C)}\Delta^{(C)}$
, the
recursive formula for the deviation
matrix $D^{(C)}$ of $Q^{(C)}$ follows from Proposition \ref{QSharpeExpress}.

\begin{prop}\label{prop_recur_dev} The deviation matrix $D^{(1)}$ of the QBD process with generator
\begin{equation}
Q^{(1)}=
\left[\begin{array}{cc}
B & A_1\\
A_{-1} & C
\end{array}\right]
\end{equation} 
 is given by 
\begin{equation}
D^{(1)}=(\boldsymbol{1\pi}^{(1)}-Q^{(1)})^{-1}-\boldsymbol{1\pi}^{(1)}.
\end{equation}
For $C\geq2$, the deviation matrix $D^{(C)}$ of the
QBD process with generator $Q^{(C)}$ is {recursively} given by
\begin{equation}
D^{(C)}=\left(\boldsymbol{1\pi}^{(C)}-I \right)
(T^{(C)})^{\#}
\left(I-E_{C-1}^{(C)}(I+\Delta^{(C)} (T^{(C)})^{\#} E_{C-1}^{(C)})^{-1}\Delta^{(C)} (T^{(C)})^{\#} \right),\label{QSharpeN+1}
\end{equation}
where $\Delta^{(C)}$ is defined in (\ref{eq:delta}),
and where $(T^{(C)})^{\#}$ is given in terms of $D^{(C-1)}$ in \eqref{eq:QsharpeTilde}.
\end{prop} 

\textit{Proof.} The group inverse of $Q^{(1)}$ comes directly from
the relationship between the group inverse of a generator and Kemeny
and Snell's fundamental matrix \cite{kemeny1976finite}, given in
Theorem 3.1 in \cite{meyer1975role}.

For any $C\geq2$, Equation (\ref{QSharpeN+1}) corresponds to Equation (\ref{QsharpeExpr2}) in Proposition \ref{QSharpeExpress} with $\tilde{D}=D^{(C)}$, $\tilde{\vc\pi}=\vc\pi^{(C)}$, $D=-(T^{(C)})^{\#}$, $K=C-1$, and $P=\Delta^{(C)}$. 
\hfill{}\cqfd

\medskip


\section{Numerical illustrations}
In this section, we compute the expected loss revenue function for various \emph{MAP/PH/1/C} queues using the results of Section \ref{mda}. We then compare the performance of the two different approaches developed in Sections \ref{mda} and \ref{pta} by computing the deviation matrix of arbitrary finite QBD processes with different numbers of phases and levels.

\subsection{The \emph{MAP/PH/1/C} queue.}

Recall that our motivation behind the computation of the expected reward function and the deviation matrices of a finite QBD process stems from our desire to compute of the expected amount of revenue lost in a \emph{MAP/PH/1/C} queue. In such a process,
\begin{itemize}
\item[\textit{(i)}] the arrival \emph{MAP} is characterised by the matrices $D_0$ and $D_1$ with $n_1$ phases;
\item[\textit{(ii)}] the service time is distributed according to a \emph{PH}$(\vc\tau,T)$ distribution of order $n_2$, with $\vc t=-T\vc1$.
\end{itemize}
The \emph{MAP/PH/1/C} system is a finite QBD process $\{\vc X(t)=(J(t),\vc\varphi(t)), t\geq 0\}$ where 
\begin{itemize}\item $0\leq J(t)\leq C$ represents the number of customers in the system at time $t$,
\item $\vc\varphi(t)=(\varphi_1(t),\varphi_2(t))$ where $0\leq\varphi_1(t)\leq n_1$ is the phase of the \emph{MAP} at time $t$, and
 $0\leq\varphi_2(t)\leq n_2$ is the phase of the \emph{PH} distribution at time $t$.
\end{itemize}
 The generator of that QBD process has a block-tridiagonal form given by \eqref{GeneratorQ}
where the block matrices are of size $n=n_1\,n_2$ and are given by
$$A_{-1}=I\otimes \vc t\cdot\vc\tau,\qquad A_0=D_0\oplus T,\qquad A_1= D_1\otimes I,$$ and $$B_0= D_0\otimes I, \qquad C_0= (D_0+D_1)\oplus T,$$where $\otimes$ and $\oplus$ denote the Kronecker product and sum, respectively.

We consider a simple example where the \emph{MAP} of arrivals corresponds to a \emph{PH} renewal process with
$$D_0=\left[\begin{array}{cc}-10 & 2\\
    1& -6\end{array}\right],\quad D_1=\left[\begin{array}{c}8 \\
    5\end{array}\right]\cdot [0.8,0.2],$$ and
the service time is \emph{PH}$(\vc\tau, T)$ with
$$\vc\tau=[0.4,0.6],\quad T=\left[\begin{array}{cc}-3 & 2\\
    1& -4\end{array}\right].$$ 
    Let $A=A_{-1}+A_0+A_1$ be the phase transition matrix associated with the QBD process, and let $\vc\alpha$ be the stationary vector of $A$. 
 Since in this example $\vc\alpha A_{-1}\vc 1< \vc\alpha A_{1}\vc 1$, we are in a \emph{high-blocking system}. This would correspond to a transient QBD process if there were no bound on the number of levels \cite{latram}. 
 
\begin{figure}[h] 
\centering
\includegraphics[angle=0, width=10cm]{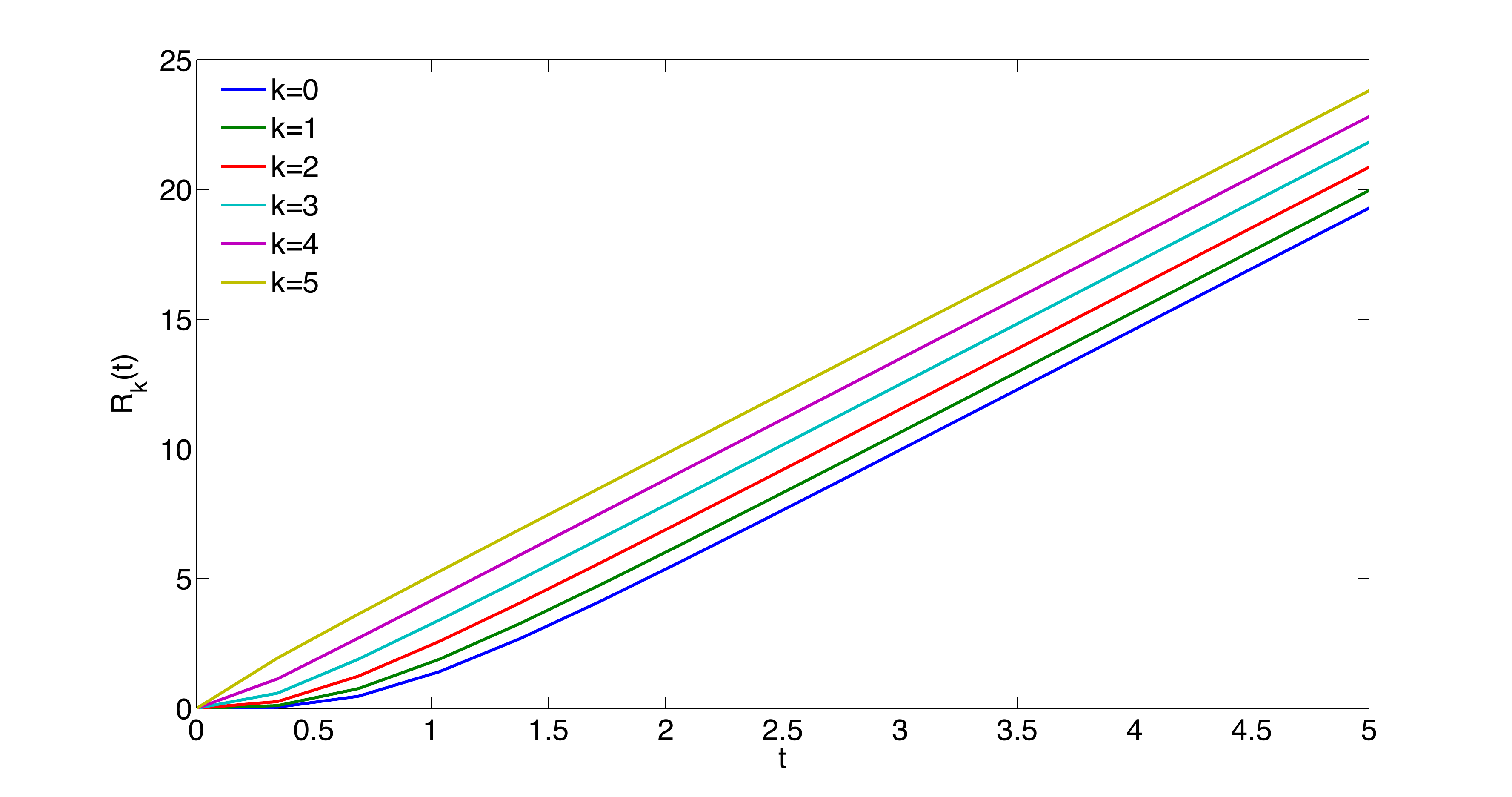} 
\includegraphics[angle=0, width=10cm]{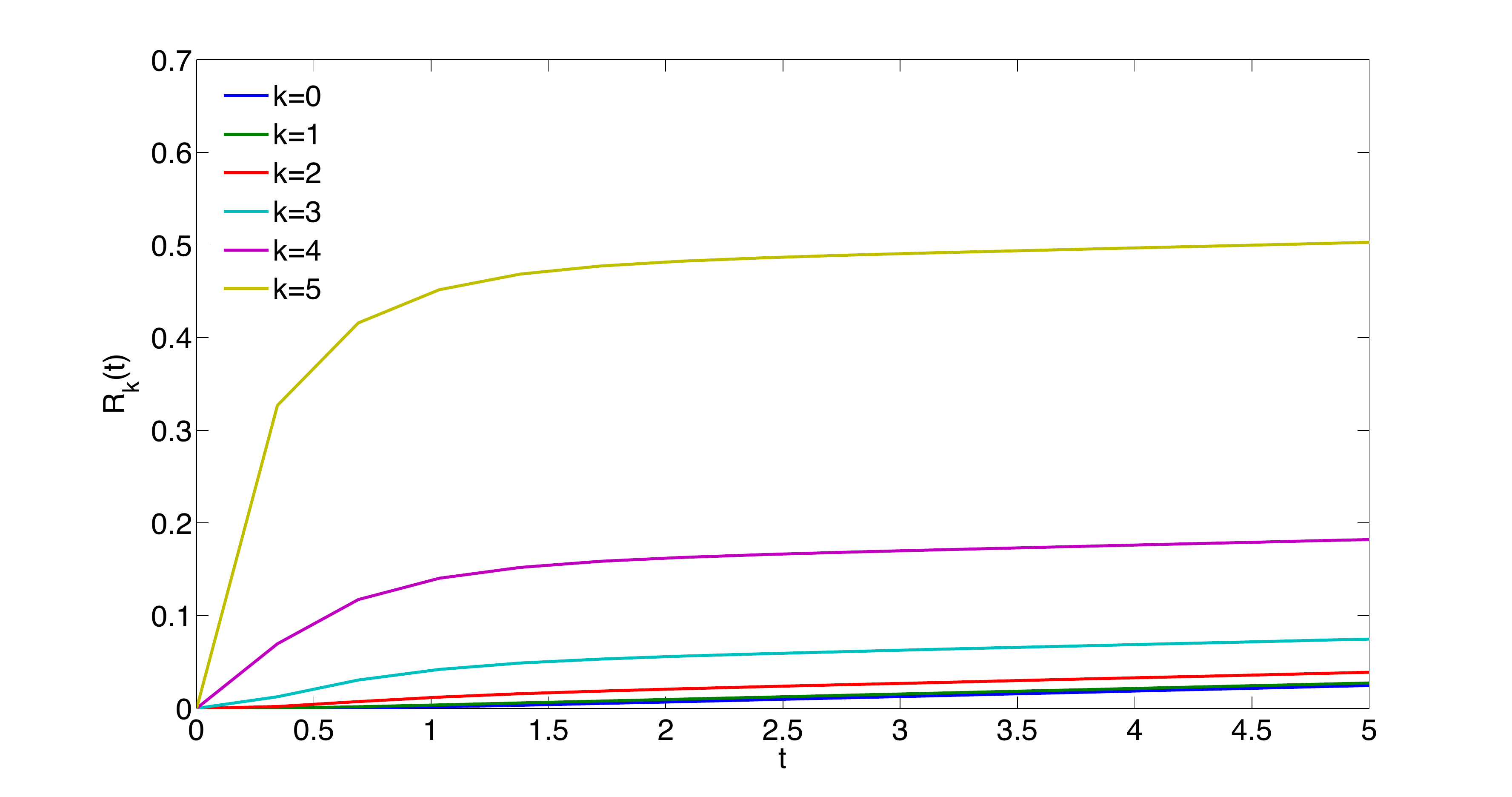}
\caption{\label{fig1}Expected lost revenue functions {corresponding to a high-blocking system (top) and a low-blocking system (bottom) with $C=5$ and $\theta=1$, for different initial queue lengths.}}
 \end{figure}

We first assume that each customer accepted in the system generates $\theta$ units of revenue. The system loses revenue when it is in level $C$ and a customer arrives, who is rejected, and does not lose any revenue otherwise. The expected reward function $\vc R(t)$ records the total expected amount of lost revenue 
over the finite time horizon $[0,t]$. In this particular case, the reward vector $\vc g$ records the \emph{loss} per unit time, per state, and we have $\vc {g}_k=\vc 0$ for $0\leq k<C$ and $\vc {g}_C=\theta A_1\vc1$. 
 Results when $C=5$ and $\theta=1$ are shown in {the upper panel of Figure \ref{fig1}} where we have assumed that the initial phase follows the distribution $\vc\alpha$. 
 
  Swapping the \emph{PH} distribution of the arrival process and service time leads to a \emph{low-blocking system}, whose results are shown in {the lower panel of Figure \ref{fig1}}. We clearly see that the expected lost revenue is much lower in the low-blocking system than in the high-blocking system, and in both cases it increases with the initial queue size, as expected.
 
\begin{figure}[h] 
\centering
\includegraphics[angle=0, width=10cm]{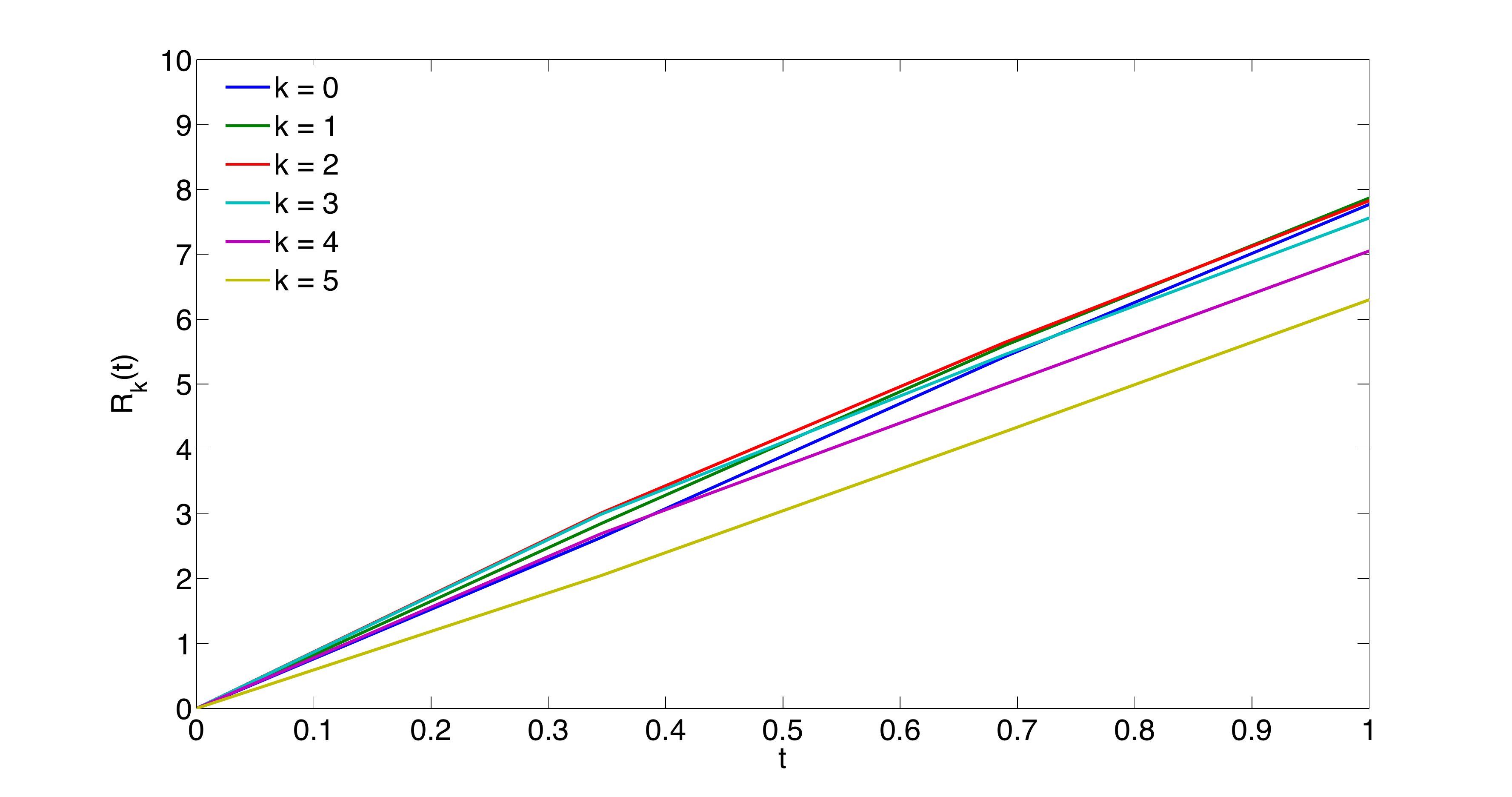} 
\includegraphics[angle=0, width=10cm]{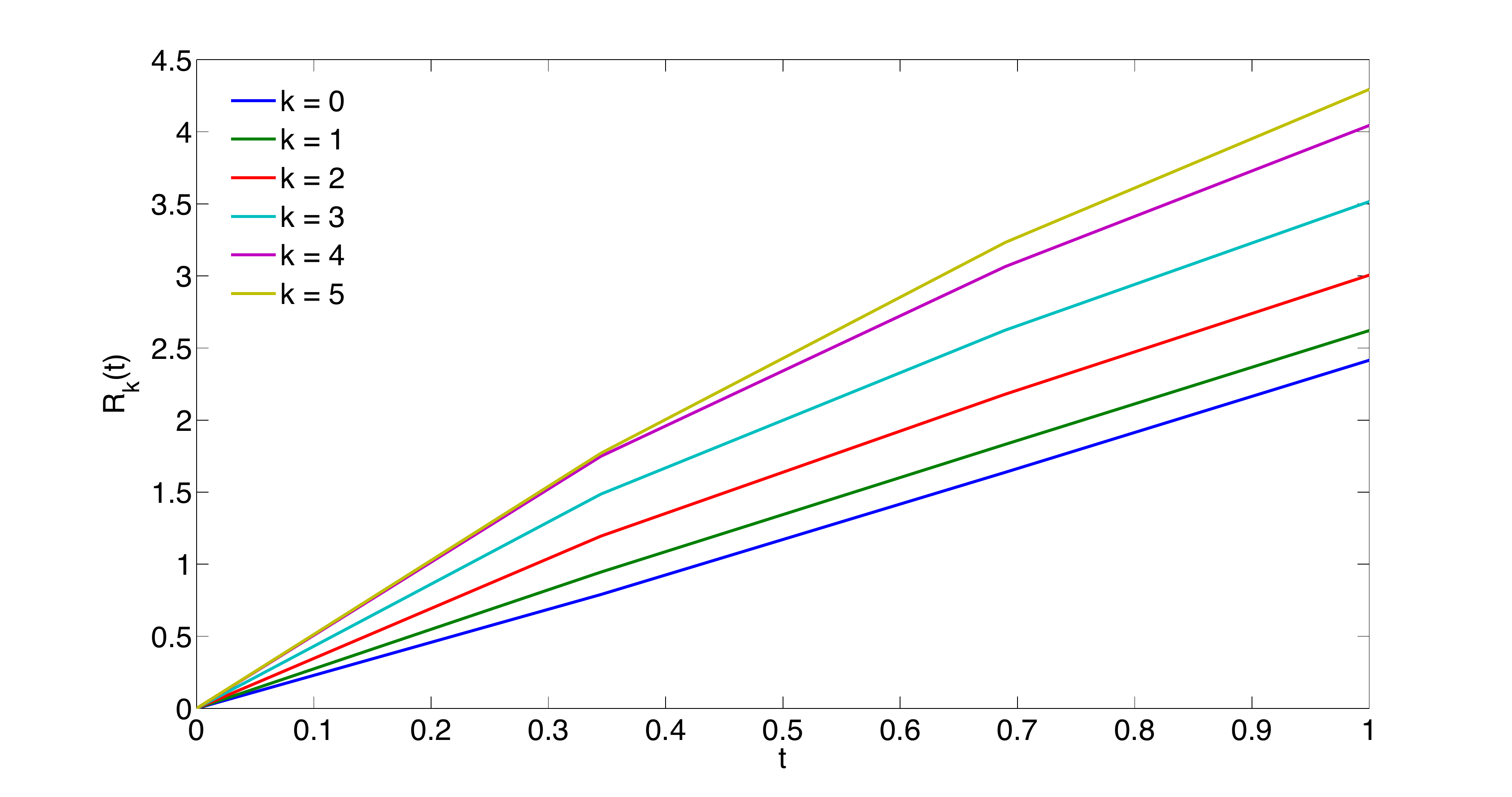}  
\caption{\label{fig3}Expected gained revenue functions {corresponding to a high-blocking system (top) and a low-blocking system (bottom) with $C=5$, $\theta=1$, and $\gamma=1$, for different initial queue lengths.}}
 \end{figure}

 We now consider a case where each customer brings a fixed amount of revenue $\gamma$ per time unit when in the system, in addition to the fixed $\theta$ units of revenue when entering the system. The reward vector then becomes $\vc g_k= \theta A_1\vc 1+\gamma k \vc 1$ for $0\leq k\leq C-1$ and $\vc g_C= \gamma C \vc 1$, and we are now interested in the expected amount of revenue \textit{gained} in $[0,t]$. The results are shown in Figure \ref{fig3} for a high-blocking system (upper panel) and a low-blocking system (lower panel).

 \subsection{Comparison of the two approaches.} 
 
 The matrix difference approach developed in Section \ref{mda} allows us to obtain the blocks of the expected reward vector and of the deviation matrices, with the advantage that the matrices involved in the expressions are at most twice the size of the phase space. However, that method requires the numerical computation of the matrices $G(s)$ and $\hat{G}(s)$ for $s\geq 0$; in addition, if we wish to compute the blocks of the deviation matrix, we first need to compute the blocks of the mean first passage time matrix, column by column.
In contrast, the perturbation theory approach developed in Section \ref{pta} allows us to obtain recursively the whole deviation matrix, but at the cost of dealing with operations on larger matrices. 
 
\begin{figure}[h] 
 \centering
\includegraphics[angle=0, width=10cm]{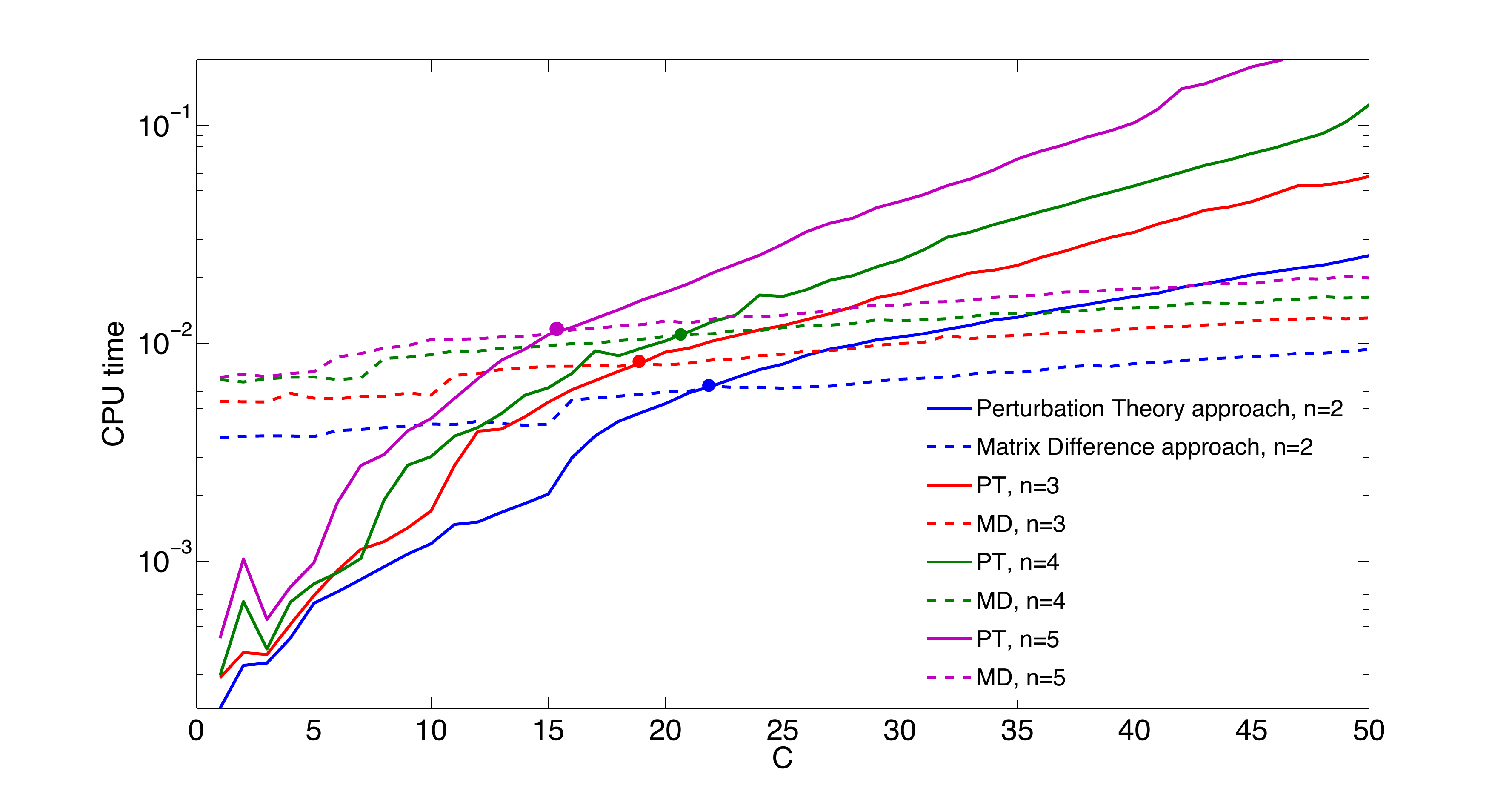}  
\caption{\label{fig5}Comparison if the CPU time corresponding to the computation of the last block-column of the deviation matrix of a QBD process with $n$ phases and maximum level $C$, using the two approaches.}
 \end{figure}

 We measured the numerical efficiency of the two approaches by comparing the CPU time corresponding to the computation of the last block-column of the deviation matrix for arbitrary QBD processes with a different number of phases and a different maximal level $C$. To that aim, we generated QBD processes with random entries, with $n=2,\ldots,5$ and $C=1,\ldots, 100$, and for each process, we ran the two methods 100 times and took the average CPU time. We repeated the same simulation several times and observed similar behaviour. The results of one experiment are shown in Figure \ref{fig5}. The CPU time increased linearly with $C$ for the matrix difference approach, while the growth is exponential for the perturbation theory approach. We see that there is always a threshold value $C^*$ such that the perturbation theoretical approach outperforms the matrix difference approach when $C<C^*$, and conversely when $C>C^*$. Note that the threshold value does not seem to show a monotonic behaviour when the number of phases $n$ increases. 
 
 In conclusion, if the objective is to compute parts of the deviation matrix, the perturbation theory approach is preferable for relatively small values of the maximal capacity, while the matrix difference approach is clearly more efficient for large values of $C$.
 
\section*{Acknowledgements} Sophie
Hautphenne and Peter Taylor would like to thank the Australian Research Council (ARC) for supporting this research through
Discovery Early Career Researcher Award DE150101044 and Laureate Fellowship FL130100039, respectively. In addition, they both acknowledge the support of the ARC Center of Excellence for Mathematical and Statistical Frontiers (ACEMS).
\\
Sarah Dendievel and Guy Latouche thank the Minist\`{e}re de la Communaut\'{e} fran\c caise de Belgique for funding this research through the ARC grant AUWB-08/13-ULB 5.
Sarah Dendievel would also like to thank the Methusalem program of the Flemish
Community of Belgium for supporting this work in part.

\end{document}